\documentclass[12pt]{amsart}
\usepackage[a4paper, left=22mm, right=22mm, top=26mm, bottom=36mm]{geometry}

\usepackage{amsmath} 					
\usepackage{amsthm} 						
\usepackage{amssymb} 					
\usepackage{amscd} 						

\usepackage{epic, eepic}
\usepackage[dvipdfm]{graphicx} 		

\usepackage[all]{xy}							

\usepackage{mathrsfs} 					
\usepackage{latexsym}						
\usepackage[russian, english]{babel}

\usepackage{enumerate}					
\usepackage{verbatim}						
\AtBeginDocument{\shorthandoff{"}}			

\newcommand{\Rational}{\mathbb{Q}}					
\newcommand{\Integer}{\mathbb{Z}}					

\newcommand{\isomarrow}{\overset{\sim}{\longrightarrow}}

\def\rank{\mathop{\mathrm{rank}}\nolimits}			
\def\corank{\mathop{\mathrm{corank}}\nolimits}		

\def\det{\mathop{\mathrm{det}}\nolimits}			
\def\Norm{\mathop{\mathrm{Nm}}\nolimits}				
\def\disc{\mathop{\mathrm{disc}}\nolimits}			

\newcommand\transpose{{}^t \!}							

\newcommand{\Sym}{{\mathrm{Sym}}}			

\def\identity{\mathop{\mathrm{id}}\nolimits}				
\def\projection{\mathop{\mathrm{pr}}\nolimits}				
\newcommand{\vertisom}{\rotatebox{270}{$\cong$}}		
\newcommand{\isomvert}{\rotatebox{90}{$\cong$}}		

\def\Image{\mathop{\mathrm{Im}}\nolimits}		
\def\Kernel{\mathop{\mathrm{Ker}}\nolimits}		
\def\Cokernel{\mathop{\mathrm{Coker}}\nolimits}		
\def\Hom{\mathop{\mathrm{Hom}}\nolimits}			
\def\Isom{\mathop{\mathrm{Isom}}\nolimits}			
\def\sheafhom{\mathop{\mathscr{H}\kern -2pt om}\nolimits}		

\def\End{\mathop{\mathrm{End}}\nolimits}			
\def\sheafend{\mathop{\mathscr{E}\kern -2pt nd}\nolimits}			
\def\Aut{\mathop{\mathrm{Aut}}\nolimits}			
\def\Isom{\mathop{\mathrm{Isom}}\nolimits}			

\def\sheafext{\mathop{\mathscr{E}\kern -2pt xt}\nolimits}			
\def\projdim{\mathop{\mathrm{pd}}\nolimits}

\def\Left{\mathop{\mathrm{L} \kern -2pt}\nolimits}				
\def\Right{\mathop{\mathrm{R} \kern -2pt}\nolimits}			

\newcommand{\Cohomology}[2]{H^{#1}\! \left( {#2} \right)}
\newcommand{\Homomorphism}[2]{\Hom_{#1}\! \left( {#2} \right)}
\newcommand{\Endomorphism}[2]{\End_{#1}\! \left( {#2} \right)}
\newcommand{\Automorphism}[2]{\Aut_{#1}\! \left( {#2} \right)}
\newcommand{\Isomorphism}[2]{\Isom_{#1}\! \left( {#2} \right)}

\newcommand{\Sheafhom}[2]{\sheafhom_{#1}\! \left( {#2} \right)}
\newcommand{\Sheafend}[2]{\sheafend_{#1}\! \left( {#2} \right)}

\newcommand{\Sheafext}[3]{\sheafext^{#1}_{#2}\! \left( {#3} \right)}

\newcommand{\strshf}{\mathcal{O}}				
\newcommand{\projsp}{\mathbb{P}}				

\def\Smooth{\mathop{\mathrm{Sm}}\nolimits}			
\def\Singular{\mathop{\mathrm{Sing}}\nolimits}			
\def\Generics{\mathop{\mathrm{Gen}}\nolimits}			

\def\codimension{\mathop{\mathrm{codim}}\nolimits}		
\def\canonical{\mathop{\mathrm{can}}\nolimits}			

\def\Supp{\mathop{\mathrm{Supp}}\nolimits}			
\def\length{\mathop{\mathrm{length}}\nolimits}
\def\depth{\mathop{\mathrm{depth}}\nolimits}		

\def\Spec{\mathop{\mathrm{Spec}}\nolimits}			
\def\sheafspec{\mathop{\mathscr{S}\kern -2pt pec}\nolimits}
\def\sheafproj{\mathop{\mathscr{P}\kern -2pt roj}\nolimits}

\def\Coh{\mathop{\mathrm{Coh}}\nolimits}		
\newcommand\GDual{\mathop{\mathrm{GD}}\nolimits}	

\def\GenLin{\mathop{\mathrm{GL}}\nolimits}			
\def\ProjGL{\mathop{\mathrm{PGL}}\nolimits}		
\def\Matrix{\mathop{\mathrm{Mat}}\nolimits}		

\makeatletter
\def\iddots{\mathinner{\mkern1mu\raise\p@
    \hbox{.}\mkern2mu\raise4\p@\hbox{.}\mkern2mu
    \raise7\p@\vbox{\kern7\p@\hbox{.}}\mkern1mu}}
\def\adots{\mathinner{\mkern2mu\raise\p@\hbox{.} 
 \mkern2mu\raise4\p@\hbox{.}\mkern1mu
 \raise7\p@\vbox{\kern7\p@\hbox{.}}\mkern1mu}}
\makeatother

\theoremstyle{definition}
\newtheorem{thm}{Theorem}[section]			
\newtheorem{prop}[thm]{Proposition}			
\newtheorem{cor}[thm]{Corollary}					
\newtheorem{lem}[thm]{Lemma}					

\newtheorem*{sketch}{Sketch of proof}

\theoremstyle{definition}
\newtheorem{dfn}[thm]{Definition}				
\newtheorem*{prf}{Proof}								
\newtheorem{rmk}[thm]{Remark}					
\newtheorem{exa}[thm]{Example}					

 \makeatletter
    
    \@addtoreset{equation}{section}
  \makeatother

\newcommand{\relmiddle}[1]{\mathrel{}\middle#1\mathrel{}}	
\newcommand{\relmid}{\relmiddle{|}}	

\newcommand{\Equref}[1]{(\ref{#1})}		

\newcommand{\CohenMacaulay}{{Cohen--Macaulay}}
\newcommand{\ArithCM}{arithmetically \CohenMacaulay}

\def\ThetaChar{\mathop{\mathrm{TC}}\nolimits}	
\def\Qdrcsp{W}
\def\NVLoci{W_{\mathrm{nv}}}			
\def\RedLoci{W_{\mathrm{gr}}}		
\def\Triples{V_{m+1, n+1}}

\title[Orbit parametrizations of theta characteristics]
{Orbit parametrizations of theta characteristics on hypersurfaces over arbitrary fields}
\author{Yasuhiro Ishitsuka}
\address{Department of Mathematics, Faculty of Science, Kyoto University, Kyoto 
606-8502, Japan}
\email{yasu-ishi@math.kyoto-u.ac.jp}
\date{\today}



\begin{document}

\maketitle

\begin{abstract}
It is well-known that theta characteristics on smooth plane curves over a field
of characteristic different from two are
in bijection with certain smooth complete intersections of three quadrics.
We generalize this bijection to possibly singular hypersurfaces of any dimension 
over arbitrary fields including those of characteristic two. 
It is accomplished in terms of linear orbits of tuples of symmetric matrices
instead of smooth complete intersections of quadrics.
{As an application of our methods, we give a description of the
projective automorphism groups of complete intersections of quadrics
generalizing Beauville's results.}  
\end{abstract}

\section{Introduction}\label{Intro}
Theta characteristics, the square roots of the canonical bundles on algebraic curves, 
are interesting objects in algebraic geometry and number theory.
They appear in several different kinds of classical problems
such as bitangents of plane quartics, determinantal representations
and Appolonius' problem (cf. \cite{Harris-Theta}, \cite[Chapter 4]{Dolgachev-CAG}).
It is well-known that, over a field of characteristic different from two, 
there is a natural bijection between 
certain smooth complete intersections of three quadrics and 
smooth plane curves with non-effective theta characteristics 
(\cite[Chapitre 6]{Beauville-Prym},
\cite[Chapter 4]{Ho-Thesis}). (In fact, Beauville also treated
the case of nodal plane curves in \cite{Beauville-Prym}.)
The purpose of this paper is to generalize this bijection.
We study more general hypersurfaces of any dimension, which may have singularities, 
over a field of arbitrary characteristic.

The main results of this paper are formulated 
in terms of linear orbits of $(m+1)$-tuples of symmetric matrices of size $n+1$
instead of complete intersections of $m+1$ quadrics in the projective space $\projsp^n.$
Of course, this formulation is equivalent to the previous one 
when $n > m$ and the characteristic of the base field is different from two. 
However, there is an essential difference when $n \le m$ or the characteristic of 
the base field is two.
For example, if $n \le m,$ the intersection of $m+1$ quadrics in $\projsp^n$ can be empty.
Hence we cannot recover the $m+1$ quadrics from their intersection. 

We fix a field $k$ of arbitrary characteristic. We fix integers $m \ge 2$ and  $n \ge 1.$ Let
\begin{align*}
	\Qdrcsp &:= 
	k^{m+1} \otimes \Sym_2 k^{n+1}\\
	&:= \left\{ M=(M_0, M_1, \dots, M_m) \; \relmid \;
	M_i \in \Matrix_{n+1}(k), \; \transpose M_i = M_i \; (i=0, 1, \dots, m) \right\}
\end{align*}
be the $k$-vector space of $(m+1)$-tuples of symmetric matrices of size $n+1$
with entries in $k.$ For an element 
\[
	M=(M_0, M_1, \dots, M_m) \in \Qdrcsp,
\]
we define its {\it discriminant polynomial} by
\[
	\disc (M) := \det (X_0 M_0 + X_1 M_1 + \dots + X_{m} M_m) 
	\in k[X_0, X_1, \dots, X_m].
\]
If $\disc(M) \neq 0,$ the discriminant polynomial $\disc(M)$ is 
a homogeneous polynomial of degree $n+1$ in $m+1$ variables $X_0, X_1, \dots, X_m.$

The $k$-vector space $\Qdrcsp$ has a natural right action
of the product of general linear groups $\GenLin_{m+1}(k) \times \GenLin_{n+1}(k).$
Concretely, for $A=(a_{i,j}) \in \GenLin_{m+1}(k), P \in \GenLin_{n+1}(k)$ and 
$M=(M_0, M_1, \dots, M_m) \in \Qdrcsp,$ we set
\[
	M \cdot (A, P) := 
	\left( \sum_{i=0}^m a_{i,0}\transpose P M_i P, \hspace{5pt}
	\sum_{i=0}^m a_{i,1}\transpose P M_i P, \hspace{5pt}\dots, \hspace{5pt}
	\sum_{i=0}^m a_{i,m}\transpose P M_i P \right).
\]
When $(A, P) \in (k^\times I_{m+1}) \times \GenLin_{n+1}(k),$
where $I_{m+1}$ is the identity matrix of size $m+1,$ 
this action preserves the discriminant polynomial of an element of $\Qdrcsp$
up to the multiplication by an element of $k^\times.$ Concretely, we have 
\[
	\disc(M \cdot (A, P)) = \det(A) \det(P)^2 \disc(M).
\]

We shall define two subsets
\[
	\RedLoci \subset \NVLoci \subset \Qdrcsp = k^{m+1} \otimes \Sym_2 k^{n+1}
\]
as follows. Let $\NVLoci$ be the subset of $(m+1)$-tuples of 
symmetric matrices whose discriminant polynomials are non-zero, and $\RedLoci$
the subset of $\NVLoci$ consisting of elements which
have no multiple factors over an algebraic closure of $k.$ 
(Here, the subscript ``nv'' stands for ``non-vanishing'', and 
the subscript ``gr'' stands for ``geometrically reduced''.)
For an element $M \in \NVLoci,$ 
the equation $(\disc (M) = 0)$ defines a hypersurface 
\[
	S \subset \projsp^m
\]
of degree $n+1$ over $k.$ 
The hypersurface $S \subset \projsp^m$ is geometrically reduced 
if and only if $M \in \RedLoci.$
The subsets $\RedLoci$ and $\NVLoci$ are stable under the 
action of $\GenLin_{m+1}(k) \times \GenLin_{n+1}(k).$

Moreover, for an element $M \in \RedLoci,$ 
we can construct an injective morphism of $\strshf_{\projsp^m}$-modules
\[
	M \colon \bigoplus_{i=0}^n \strshf_{\projsp^m}(-2) 
	\longrightarrow \bigoplus_{i=0}^n \strshf_{\projsp^m}(-1)
\]
and a coherent $\strshf_S$-module $\mathcal{M}$ by the following short exact sequence
of $\strshf_{\projsp^m}$-modules:
\[
	\xymatrix{
		0 \ar[r] & {\displaystyle \bigoplus_{i=0}^{n} \strshf_{\projsp^m}(-2)} \ar[r]^{M}
		& {\displaystyle \bigoplus_{i=0}^{n} \strshf_{\projsp^m}(-1)} \ar[r]^(.7){p}
		& \mathcal{M} \ar[r] & 0.
	}
\]
The coherent $\strshf_S$-module $\mathcal{M}$ as above is called 
a {\it non-effective theta characteristic} on $S.$ 
(It is {\ArithCM} and pure of dimension $m-1.$
The length of it at each generic point of $S$ is equal to one.
It also satisfies $\Cohomology{0}{S, \mathcal{M}} = 0,$ 
and has a certain duality quasi-isomorphism $\lambda.$
See Definition \ref{TC} for details.)
 
Let $\ThetaChar_{m+1, n+1}(k)$ be the set of equivalence classes of 
pairs $(S, \mathcal{M})$ 
which consist of a geometrically reduced hypersurface 
$S \subset \projsp^m$ of degree $n+1$ over $k$ and a non-effective theta characteristic
 $\mathcal{M}$ on $S.$ 
Here, two pairs $(S, \mathcal{M}), (S',\mathcal{M}')$ are said to be
{\it equivalent} if $S = S'$ and 
$\mathcal{M}, \mathcal{M}'$ are isomorphic 
as $\strshf_S$-modules (see Subsection \ref{main} for details).
When a pair $(S, \mathcal{M})$ comes from an element $M \in \RedLoci,$
the equivalence class $[(S, \mathcal{M})]$ of a pair $(S, \mathcal{M})$ is 
determined by the $(k^\times I_{m+1}) \times \GenLin_{n+1}(k)$-orbit of $M.$ 
We put 
\[
	\Phi_{m+1, n+1}(M) := [(S, \mathcal{M})].
\]
Hence we have obtained a map 
\[
	\Phi_{m+1, n+1} \colon \RedLoci
	/((k^\times I_{m+1}) \times \GenLin_{n+1}(k)) \to \ThetaChar_{m+1, n+1}(k).
\]

In this paper, we prove the surjectivity of the map $\Phi_{m+1, n+1}$
and study the structure of the fibers of $\Phi_{m+1, n+1}.$

\begin{thm}\label{main-theorem}
	Let $k$ be a field of arbitrary characteristic.
	Let $m, n$ be integers satisfying $m \ge 2$ and $n \ge 1.$
	\begin{enumerate}
		\item The map $\Phi_{m+1, n+1}$ is surjective.
		\item For an element $[(S, \mathcal{M})] \in \ThetaChar_{m+1, n+1}(k),$
		there exists a finite dimensional commutative {\'e}tale $k$-algebra $L$
		such that the fiber $\Phi_{m+1, n+1}^{-1}\left( [( S, \mathcal{M} )] \right)$ has 
		a simply transitive action of the group $L^\times / k^\times L^{\times 2}.$ 
	\end{enumerate}
\end{thm}

In the second statement of Theorem \ref{main-theorem},
the group $L^\times / k^\times L^{\times 2}$ is 
the quotient of the multiplicative group $L^\times$ by the subgroup
\[
	k^\times L^{\times 2} := \left\{ ab^2 \relmid 
	a \in k^\times, b \in L^\times \right\}. 
\]
We also obtain some results on tuples of symmetric matrices 
defining hypersurfaces which are not necessarily geometrically reduced
(see Corollary \ref{Bij-cor1},
Corollary \ref{Bij-cor2} and Proposition \ref{fiber}). 

From Theorem \ref{main-theorem}, we can recover some results of Beauville and Ho
(cf.\ \cite{Beauville-Prym}, \cite{Ho-Thesis}):

\begin{cor}[{See Corollary \ref{BeauTheo}}]\label{Beauville-theorem}
	The fiber $\Phi_{m+1, n+1}^{-1}([(S, \mathcal{M})])$ is a singleton 
	if {\it at least one} of the following conditions is satisfied:
   \begin{itemize}
  		\item the base field $k$ is separably closed of characteristic different from two, or
		\item the base field $k$ is perfect of characteristic two, or
		\item the hypersurface $S \subset \projsp^m$ is geometrically integral.
   \end{itemize}
\end{cor}

In fact, it is easy to see that the group $L^\times / k^\times L^{\times 2}$ is trivial
if at least one of the above conditions is satisfied.

Beauville proved Corollary \ref{Beauville-theorem} 
when $k$ is algebraically closed of characteristic 
different from two, $m=2, n \ge 3,$ $S \subset \projsp^2$ is a nodal plane curve
and one of the elements $M \in \Phi^{-1}_{m+1, n+1}([(S, \mathcal{M})])$ 
in the fiber defines a smooth complete intersection of
three quadrics in $\projsp^n$ ({\cite[Proposition 6.19]{Beauville-Prym}}).
Also, Corollary \ref{Beauville-theorem} is proved by Ho when $m=2, n \ge 2,$
$k$ is a field where $3n(n-1)$ is invertible
and $S \subset \projsp^2$ is a smooth plane curve 
({\cite[Theorem 4.12]{Ho-Thesis}}). 
Note that, if the set of conditions imposed by Beauville or Ho is satisfied,
the fiber of $\Phi_{m+1, n+1}$ is a singleton.
However, the map $\Phi_{m+1, n+1}$ is not injective in general.
Actually, a fiber can have infinitely many elements.
(See Example \ref{infinite}.)

Another result of this paper concerns the projective automorphism groups
of complete intersections of quadrics.
For this application, we assume that the characteristic of $k$ is different from two and
$n>m\ge 2.$
Let $(S, \mathcal{M})$ be a pair of 
a geometrically reduced hypersurface $S \subset \projsp^m$ and 
a non-effective theta characteristic $\mathcal{M}$ on $S.$
We take an element $M \in \Phi_{m+1,n+1}^{-1}([(S, \mathcal{M})]).$
Define the closed subvariety $X_Q$ of $\projsp^n$ by
\[
	X_Q := \left\{ x \in \projsp^n \relmid \transpose x M_i x = 0 
	\; (i=0,1, \dots, m) \right\}.
\]
Here, we identify points on $\projsp^n$ and
$(n+1)$-dimensional column vectors.
The projective isomorphism class of $X_Q$ depends only on the 
$\GenLin_{m+1}(k) \times \GenLin_{n+1}(k)$-orbit of $M.$
The choice of $M$ defines a duality quasi-isomorphism $\lambda$ of $\mathcal{M}$
and an {\'e}tale $k$-algebra $L$ in Theorem \ref{main-theorem}.
Let $G$ be the kernel of the homomorphism
\[
	L^\times / k^\times \longrightarrow L^\times/k^\times \quad ; \quad
	x \mapsto x^2.
\]
Define
\(
	\Automorphism{\projsp^m}{S, \mathcal{M}, \lambda}
\)
to be the subgroup of the projective automorphism group $\Automorphism{\projsp^m}{S}$ 
of $S$ consisting of elements $\nu \in \Automorphism{\projsp^m}{S}$ 
fixing the equivalence class of $(\mathcal{M}, \lambda)$
(see Section \ref{ProjAut} for details). 
Then we have the following short exact sequence.

\begin{thm}[See Theorem \ref{ProjAutComm}]\label{ProjAutComm0}
	Assume that the characteristic of $k$ is different from two and $n>m\ge 2.$
	Let $M \in \Phi_{m+1,n+1}^{-1}([(S, \mathcal{M})])$ 
	be an element defining a complete intersection $X_Q$ of $m+1$ quadrics in $\projsp^n.$
	Then there exists a short exact sequence of the following form:
		\[
			\xymatrix{
				0 \ar[r] &
				G \ar[r] &
				\Automorphism{\projsp^n}{X_Q} \ar[r] &
				\Automorphism{\projsp^m}{S, \mathcal{M}, \lambda} \ar[r] &
				0.
			}
		\]
\end{thm}
Beauville proved Theorem \ref{ProjAutComm0} when $m=2, n \ge 3,$ 
$k$ is algebraically closed and $X_Q$ is a smooth complete intersection
of three quadrics (\cite[Proposition 6.19]{Beauville-Prym}).

Let us give a sketch of the proof of Theorem \ref{main-theorem}.
First, we prove a rigidified version of a result of Beauville 
(\cite[Proposition 1.11]{Beauville-Det}).
Precisely, after fixing a quasi-isomorphism of complexes of 
coherent $\strshf_{\projsp^m}$-modules 
\[
	c \colon \omega_{\projsp^m} \isomarrow \strshf_{\projsp^m}(-m-1)[m],
\] 
where $\omega_{\projsp^m}$ is the dualizing complex on $\projsp^m,$
we construct a bijection between the set of $(m+1)$-tuples of 
symmetric matrices (not equivalence classes of them!)\ 
and the set of equivalence classes of 
certain coherent $\strshf_{\projsp^m}$-modules with rigidification data
(see Theorem \ref{bijection1}).
It naturally induces a bijection between $\RedLoci$ and a set
which has a natural surjection onto $\ThetaChar_{m+1, n+1}(k).$
Then we prove Theorem \ref{main-theorem} by studying the fiber of the surjection.

The present paper is organized as follows.
We recall basic results on minimal resolutions of 
coherent $\strshf_{\projsp^m}$-modules in Section \ref{Preliminaries}.
In Section \ref{bijection}, we give a rigidified bijection as explained above 
in a general situation by using Grothendieck duality and 
minimal resolutions of coherent $\strshf_{\projsp^m}$-modules.
In Section \ref{GenLin},
we study the actions of $\GenLin_{m+1}(k) \times \GenLin_{n+1}(k)$
on $\RedLoci$ and $\ThetaChar_{m+1, n+1}(k).$ 
In Section \ref{fibers}, we study symmetric quasi-isomorphisms
via the methods of Section \ref{GenLin}.
In Section \ref{Curves}, we restrict our results to the case of 
geometrically reduced hypersurfaces
using the notion of theta characteristics on hypersurfaces
following Mumford, Harris and Piontkowski 
(\cite{Mumford-Theta}, \cite{Harris-Theta}, \cite{Piontkowski-Theta}).
Then we prove Theorem \ref{main-theorem}.
Section \ref{ProjAut} is devoted to an application of our methods to
the projective automorphism groups of complete intersections of quadrics
in terms of the endomorphism ring of theta characteristics.

\subsection*{Acknowledgements}
The author would like to thank Tetsushi Ito for copious comments and continual encouragement.
This work is supported by JSPS KAKENHI Grant Number 13J01450.

\subsection*{Notation}
We work over a field $k$ of arbitrary characteristic 
except in Section \ref{ProjAut}. In Section \ref{ProjAut},
we assume the characteristic of $k$ is different from two. 
The $k$-vector space of symmetric matrices of size $n+1$ with entries in $k$
is denoted by $\Sym_2 k^{n+1}.$
Hence an element of $k^{m+1} \otimes \Sym_2 k^{n+1}$ is identified with
an $(m+1)$-tuple of symmetric matrices of size $n+1$ with entries in $k.$
For a scheme $X$ over $k$ of finite type,
we write the set of singular points, smooth points and generic points on $X$ 
as $\Singular(X), \Smooth(X)$ and $\Generics(X),$ respectively.
For a point $p \in X,$ the local ring at $p$ is denoted by $\strshf_{X, p}.$
We use the symbol $\omega_X$ to denote the {\it dualizing complex} on $X$
rather than the dualizing sheaf (cf.\ \cite{Hartshorne-Residue}).
For an object $\mathcal{F} \in D(\Coh(\strshf_X))$ in the derived category of 
complexes of coherent $\strshf_X$-modules,  
let $\mathcal{F}[n]$ denote the degree $n$ shift defined by 
\(
	(\mathcal{F}[n])_i := \mathcal{F}_{n+i}.
\)
For a morphism $h$ between complexes of $\strshf_{X}$-modules,
we denote the induced morphism between cohomologies by the same symbol $h.$

\section{Preliminaries on minimal resolutions of sheaves on $\projsp^m$}\label{Preliminaries}
We recall some notions and basic results on minimal resolutions of 
coherent $\strshf_{\projsp^m}$-modules. 
We also recall the notion and properties of {\ArithCM} $\strshf_{\projsp^m}$-modules
(cf.\ \cite[Subsection 1.5]{Beauville-Det}, \cite{Eisenbud-Syzygies}).

We fix a field $k$ of arbitrary characteristic and 
integers $m \ge 2$ and $n \ge 1.$ Let 
\[
	R := k[X_0,X_1,\dots, X_m] 
\]
be the polynomial algebra over $k$ in $m+1$ variables $X_0, X_1, \dots, X_m.$
The ring $R$ is a graded $k$-algebra, and
\[
	\mathfrak{m}_R := (X_0, X_1, \dots, X_m)R
\]
is a graded maximal ideal of $R.$

A {\it free resolution} of a finitely generated 
graded $R$-module $N$ is an exact sequence of the form
\[
	\xymatrix{
		\cdots \ar[r]^{\delta_2} & F_1 \ar[r]^{\delta_1}& 
		F_0 \ar[r]^(.37){\delta_0} & F_{-1} :=N \ar[r] & 0,
	}
\]
where $F_i$ is a graded free $R$-module 
and $\delta_i$ is an $R$-homomorphism preserving degree for each $i \ge 0.$
It is {\it minimal} if the image of $\delta_i$ is contained in the $R$-submodule 
$\mathfrak{m}_R F_{i-1}$ of $F_{i-1}$ for each $i \ge 1.$
For an integer $e \in \Integer,$ we denote by $R(e)$ the graded free $R$-module of rank one 
generated by an element of degree $-e.$
The minimal free resolution of $N$ is said to be {\it pure} if
each graded piece $F_i$ ($i \ge 0$) is generated by elements of the same degree;
in other words, there exist integers $e_i$ and $n_i$
such that $F_i \cong R(e_i)^{n_i}$ for each $i \ge 0.$

A {\it graded locally free resolution} of a coherent $\strshf_{\projsp^m}$-module
$\mathcal{N}$ is an exact sequence of the form
\begin{equation}\label{Resol}
	\xymatrix{
		\cdots \ar[r]^{\delta_2} & \mathcal{F}_1 \ar[r]^{\delta_1}& 
		\mathcal{F}_0 \ar[r]^(.37){\delta_0} & \mathcal{F}_{-1} := \mathcal{N} \ar[r] & 0,
	}
\end{equation}
where $\mathcal{F}_i$ is a direct sum of line bundles on $\projsp^m$
for each $i \ge 0.$
It is {\it minimal} if the image of the sequence \Equref{Resol} by the functor $\Gamma_*$  is
a minimal free resolution of 
a graded $\Gamma_*(\strshf_{\projsp^m}) \cong R$-module, where 
\[
	\Gamma_*(\mathcal{F}) := 
	\bigoplus_{j \in \Integer} \Cohomology{0}{\projsp^m, \mathcal{F}(j)}
\]
for an $\strshf_{\projsp^m}$-module $\mathcal{F}.$
The minimal graded locally free resolution of $\mathcal{F}$ is {\it pure} if
the image of the sequence \Equref{Resol} by the functor $\Gamma_*$ is
a pure minimal free resolution of $\Gamma_*(\mathcal{F}).$
Note that $\Gamma_*(\strshf_{\projsp^m}(e))$ is isomorphic to $R(e)$
as a graded $\Gamma_*(\strshf_{\projsp^m}) \cong R$-module.

\begin{lem}
	A minimal free resolution of a finitely generated graded $R$-module $N$ 
	exists and it is unique up to isomorphism in the category of graded $R$-modules.
	Similarly, if the set of the associated points of coherent $\strshf_{\projsp^m}$-module 
	$\mathcal{N}$ does not contain any closed point in $\projsp^m,$
	a minimal graded locally free resolution of $\mathcal{N}$ exists and 
	unique in the category of coherent $\strshf_{\projsp^m}$-modules.
\end{lem}
\begin{prf}
	For a graded $R$-module $N,$ see \cite[Exercise 20.1]{Eisenbud-Rings}.
	If a coherent $\strshf_{\projsp^m}$-module $\mathcal{N}$ has no associated points
	of dimension zero, the graded $R$-module $\Gamma_*(\mathcal{N})$ is finitely generated.
	Thus the assertion follows from the case of graded modules. \qed
\end{prf}
\begin{lem}\label{uniqueness}
	If a minimal free resolution of a finitely generated graded $R$-module $N$ is pure, 
	any endomorphism of $N$ is uniquely lifted to
	an endomorphism of each graded piece of the minimal resolution.
	If the set of the associated points of coherent $\strshf_{\projsp^m}$-module 
	$\mathcal{N}$ does not contain any closed point in $\projsp^m,$
	the same conclusion holds for endomorphisms of $\mathcal{N}.$
\end{lem}
\begin{prf}
	We only prove the case of graded $R$-modules because
	the case of coherent $\strshf_{\projsp^m}$-modules can be proved in a similar way.
	By the purity of minimal resolution, the graded $R$-module $N/{\mathfrak{m}_R}N$
	concentrates in a single degree $d.$
	Then we have isomorphisms of $k$-vector spaces
	\[
		N_d \cong (N/\mathfrak{m}_RN)_d \cong (F_0/\mathfrak{m}_RF_0)_d
		\cong (F_0)_d.
	\]
	The endomorphism of $N$ preserving degree induces an endomorphism of $N_d.$
	Hence it induces an endomorphism of $(F_0)_d.$
	Since $F_0$ is pure, $F_0$ is generated by $(F_0)_d$ as an $R$-module
	and any endomorphism of $(F_0)_d$ is uniquely lifted to $F_0.$
	By repeating this arguments, we have a unique lift to an endomorphism of each graded piece.
\qed\end{prf}

\begin{dfn}[{\cite[Subsection 1.1]{Beauville-Det}}]
A coherent $\strshf_{\projsp^m}$-module $\mathcal{F}$ is {\it \ArithCM} 
if it satisfies the following two conditions:
\begin{itemize}
	\item $\mathcal{F}$ is a {\CohenMacaulay} $\strshf_{\projsp^m}$-module, 
	that is, $\mathcal{F}_x$ is a {\CohenMacaulay} $\strshf_{\projsp^m, x}$-module
	for each point $x \in \projsp^m,$ and
	\item 
	\(
		\Cohomology{i}{\projsp^m, \mathcal{F}(j)} = 0
	\)
	for any $j \in \Integer$ and $1 \le i \le \dim \Supp (\mathcal{F}) -1.$ 
\end{itemize}
\end{dfn}

\begin{dfn}[{\cite[Definition 1.1.2]{Huybrechts-Lehn-Moduli}}]
A non-zero coherent $\strshf_{\projsp^m}$-module $\mathcal{F}$ is
{\it pure of dimension $d$} if, for any non-zero $\strshf_{\projsp^m}$-submodule 
$\mathcal{G} \subset \mathcal{F},$ the dimension of the support $\Supp(\mathcal{G})$ is
equal to $d.$
\end{dfn}
We quote the following characterization of pure sheaves of dimension $d.$
Note that there is a degree $m$ shift in the statement
because $\omega_{\projsp^m}$ denotes the {\it dualizing complex} on $\projsp^m$ 
in this paper whereas it denotes the {\it dualizing sheaf} on $\projsp^m$
in \cite{Huybrechts-Lehn-Moduli}.

\begin{prop}[{\cite[Proposition 1.1.10]{Huybrechts-Lehn-Moduli}}]\label{purity}
	Let $\mathcal{F}$ be a non-zero coherent $\strshf_{\projsp^m}$-module
	with support $\Supp(\mathcal{F})$ of dimension $d.$
	Then the following two conditions are equivalent:
	\begin{itemize}
		\item $\mathcal{F}$ is pure of dimension $d.$ 
		\item \(
		\codimension(\Supp(\Sheafext{q-m}{\projsp^m}
		{\mathcal{F}, \omega_{\projsp^m}})) \ge q+1
		\)
		for any $q > m-d.$ 
	\end{itemize}
\end{prop}

We introduce the following proposition which is essentially equivalent to 
\cite[Proposition 1.11]{Beauville-Det}.

\begin{prop}[{\cite[Proposition 1.11]{Beauville-Det}}]\label{ACM-and-resolution}
	For a coherent $\strshf_{\projsp^m}$-module $\mathcal{F},$
	the following two conditions on $\mathcal{F}$ are equivalent:
	\begin{itemize}
		\item\label{ACM1} The coherent $\strshf_{\projsp^m}$-module 
		$\mathcal{F}$ has a minimal graded locally free resolution of the following form
		for some $r \ge 0$:
			\begin{align}\label{ACME}
				\xymatrix{
					0 \ar[r] & {\displaystyle \bigoplus_{i=0}^r \strshf_{\projsp^m}(-2)} \ar[r]^M
					& {\displaystyle \bigoplus_{i=0}^r \strshf_{\projsp^m}(-1)} \ar[r]
					& \mathcal{F} \ar[r] & 0.
				}
			\end{align}
	
		\item\label{ACM2} The coherent $\strshf_{\projsp^m}$-module $\mathcal{F}$ is
			{\ArithCM} and pure of dimension $m-1.$ It also satisfies
			\[
					\Cohomology{0}{\projsp^m, \mathcal{F}} =
					\Cohomology{m-1}{\projsp^m, \mathcal{F}(2-m)} = 0.
				\]
	\end{itemize}
	The map $M$ as in \Equref{ACME} can be described as
	a square matrix of size $r+1$ with entries in 
	$\Cohomology{0}{\projsp^m, \strshf_{\projsp^m}(1)}.$	
	If $\mathcal{F}$ satisfies these equivalent conditions, 
	the support $\Supp(\mathcal{F}) \subset \projsp^m$ of $\mathcal{F}$
	is a hypersurface defined by the equation $(\det(M)=0).$
\end{prop}

\begin{sketch}
We give a brief sketch of the proof of Proposition \ref{ACM-and-resolution} 
because we need a rigidified version 
of this proposition in Section \ref{bijection}.
If $\mathcal{F}$ satisfies the first condition, $\mathcal{F}$ is 
{\ArithCM} and pure of dimension $m-1$ 
by \cite[Proposition 1.2]{Beauville-Det}.
By the long exact sequence of cohomology, we have
\[
	\Cohomology{0}{\projsp^m, \mathcal{F}} = 
	\Cohomology{m-1}{\projsp^m, \mathcal{F}(2-m)} = 0.
\]
Moreover, by the long exact sequence of cohomology, 
we have an exact sequence
\[
	\xymatrix{
		{\displaystyle \bigoplus_{i=0}^r \Sheafext{q-m-1}{\projsp^m}
		{\strshf_{\projsp^m}(-2), \omega_{\projsp^m}}} \ar[r] &
		\Sheafext{q-m}{\projsp^m}
		{\mathcal{F}, \omega_{\projsp^m}} \ar[r] &
		{\displaystyle \bigoplus_{i=0}^r \Sheafext{q-m}{\projsp^m}
		{\strshf_{\projsp^m}(-1), \omega_{\projsp^m}}},
	}
\]
where $\omega_{\projsp^m}$ denotes the dualizing complex on $\projsp^m.$
Recall that $\omega_{\projsp^m} \cong \strshf_{\projsp^m}(-m-1)[m].$
Since
\[
	\Sheafext{q-m-1}{\projsp^m}{\strshf_{\projsp^m}(-2), \omega_{\projsp^m}} \cong 
	\Sheafext{q-1}{\projsp^m}{\strshf_{\projsp^m},\strshf_{\projsp^m}(-m+1)}
	= 0
\]
for $q>1$ and
\[
	\Sheafext{q-m}{\projsp^m}{\strshf_{\projsp^m}(-1), \omega_{\projsp^m}} \cong 
	\Sheafext{q}{\projsp^m}{\strshf_{\projsp^m},\strshf_{\projsp^m}(-m)}
	= 0
\]
for $q\ge 1$ (\cite[III, Proposition 6.3]{Hartshorne-AG}), we have
\[
	\Sheafext{q-m}{\projsp^m}{\mathcal{F}, \omega_{\projsp^m}} = 0
\]
for any $q > 1.$
Thus we see that $\mathcal{F}$ is pure of dimension $m-1$ by Proposition \ref{purity}.

Conversely, assume that $\mathcal{F}$ satisfies the second condition. 
Since the coherent $\strshf_{\projsp^m}$-module $\mathcal{F}$ is arithmetically
Cohen--Macaulay, $\Gamma_*(\mathcal{F})$ is a Cohen--Macaulay $R$-module
by \cite[Proposition 1.2]{Beauville-Det}. 
For a non-zero finitely generated graded $R$-module $N,$ 
the Auslander--Buchsbaum formula
for the polynomial algebra $R$ states that
\[
	\projdim(N) + \depth(\mathfrak{m}_R, N) = \depth (\mathfrak{m}_R, R) = m+1,
\] 
where $\projdim(N)$ (resp.\ $\depth(\mathfrak{m}_R, N)$) is the projective dimension
(resp.\ $\mathfrak{m}_R$-depth) of $N$ (\cite[Exercise 19.8]{Eisenbud-Rings}).
Since $\mathcal{F}$ is {\ArithCM} and pure of dimension $m-1,$
the graded $R$-module $\Gamma_*(\mathcal{F})$ is {\CohenMacaulay} of depth $m.$
Hence, the projective dimension of $\Gamma_*(\mathcal{F})$
is one. There exists a minimal free resolution of the following form for some $r \ge 0:$
\[
	\xymatrix{
		0 \ar[r] & {\displaystyle \bigoplus_{i =0}^r R(e_i)} \ar[r] &
		{\displaystyle \bigoplus_{i =0}^r R(d_i)} \ar[r] & \Gamma_*(\mathcal{F}) 
		\ar[r] & 0.
	}
\]
Therefore, the coherent $\strshf_{\projsp^m}$-module $\mathcal{F}$ has 
a minimal graded locally free resolution of the following form:
\[
	\xymatrix{
		0 \ar[r] & {\displaystyle \bigoplus_{i =0}^r \strshf_{\projsp^m}(e_i)} \ar[r] &
		{\displaystyle \bigoplus_{i =0}^r \strshf_{\projsp^m}(d_i)} \ar[r] 
		& \mathcal{F} 
		\ar[r] & 0.
	}
\]
Since $\Cohomology{p}{\projsp^m, \strshf_{\projsp^m}(q)} = 0$ for all $1 \le p \le m-1$
and all $q \in \Integer$ (\cite[III, Theorem 5.1]{Hartshorne-AG}),
by the long exact sequence of cohomology,
we have $\Cohomology{p}{\projsp^m, \mathcal{F}(q)} = 0$
for all $1 \le p \le m-2$ and all $q \in \Integer.$
Since $\dim \Supp(\mathcal{F}) = m-1,$
we have $\Cohomology{p}{\projsp^m, \mathcal{F}(q)} = 0$
for all $p \ge m$ and all $q \in \Integer.$
By assumption, we also have $\Cohomology{m-1}{\projsp^m, \mathcal{F}(2-m)} = 0.$

Therefore, we have $\Cohomology{p}{\projsp^m, \mathcal{F}(1-p)}=0$
for all $p \ge 1,$ and $\mathcal{F}$ is {\it 1-regular} in the sense of Mumford
(\cite[Lecture 14]{Mumford-Lectures}).
Hence $\mathcal{F}(1)$ is generated by its global sections,
and the map
\[
	\Cohomology{0}{\projsp^m, \mathcal{F}(q)} \otimes
	\Cohomology{0}{\projsp^m, \strshf_{\projsp^m}(1)} \longrightarrow
	\Cohomology{0}{\projsp^m, \mathcal{F}(q+1)} \quad (q \ge 1)
\]
is surjective (\cite[Corollary 4.18]{Eisenbud-Syzygies}, 
\cite[Lecture 14]{Mumford-Lectures}).
Since $\Cohomology{0}{\projsp^m, \mathcal{F}}=0$ by assumption,
the graded $R$-module $\Gamma_*(\mathcal{F})$ is generated by elements of degree one.
By the minimality of resolution, we have $d_i = -1$ and $e_i \le -2$ for all $i.$
By the long exact sequence of cohomology,
we have $\Cohomology{m}{\projsp^m, \strshf_{\projsp^m}(e_i+2-m)}=0$ for all $i.$
Hence we have $e_i + 2-m \ge -m.$
(Recall that 
$\Cohomology{m}{\projsp^m, \strshf_{\projsp^m}(q)} = 0$ if and only if $q \ge -m$
(\cite[III, Theorem 5.1]{Hartshorne-AG}).)
We conclude that $e_i = -2$ for all $i,$
and $\mathcal{F}$ has a graded locally free resolution of the desired form.
\qed \end{sketch}

\section{A bijection between tuples of symmetric matrices and sheaves on $\projsp^m$}\label{bijection}
In this section, 
we establish a bijection related to Theorem \ref{main-theorem}.
It treats coherent $\strshf_{\projsp^m}$-modules.
It is a rigidified version of Beauville's results 
(\cite[Proposition 1.11]{Beauville-Det}).
The proofs in this section are similar to those in \cite{Beauville-Det} and
\cite[Chapter 4]{Dolgachev-CAG}.

As in the previous section, 
we fix a field $k$ of arbitrary characteristic and integers $m \ge 2$ and $n \ge 1.$
We fix a quasi-isomorphism of complexes of 
coherent $\strshf_{\projsp^m}$-modules
\[
	c \colon \omega_{\projsp^m} \isomarrow \strshf_{\projsp^m}(-m-1)[m].
\]

\subsection{The statement of the first bijection}\label{bij}
Let us introduce some notation on morphisms between complexes
of coherent sheaves in the derived category.

Let $X$ be a scheme over $k$ of finite type.
For a morphism $h \colon \mathcal{F} \to \mathcal{G}$ of bounded complexes of 
coherent $\strshf_{X}$-modules,
the {\it transpose morphism} of $h$ is defined by 
\[
	\transpose h \colon \Right\Sheafhom{X}{\mathcal{G}, \omega_X} \longrightarrow
	\Right\Sheafhom{X}{\mathcal{F}, \omega_X} \quad ; \quad 
	g \mapsto g \circ h,
\]
where $\omega_X$ is the {\it dualizing complex} on $X.$ 
If $X$ is smooth over $k$ of dimension $\dim X,$
we have a canonical quasi-isomorphism
\(
	\omega_X \cong \Omega_{X/k}^{\dim X}[\dim X],
\)
where $\Omega_{X/k}^{\dim X}$ is the canonical sheaf on $X$
(\cite[Theorem 4.1]{Hartshorne-Residue}).
The following canonical homomorphism
\[
	\canonical_{\mathcal{F}, \mathcal{G}} \colon 
	\mathcal{F} \longrightarrow \Right \Sheafhom{X}
	{\Right \Sheafhom{X}{\mathcal{F}, \mathcal{G}}, \mathcal{G}} \quad ; \quad
	g \mapsto (s \mapsto s(g))
\]
is also denoted by ``$\canonical$'' if there is no danger of confusion.
Note that $\canonical_{\mathcal{F}, \omega_X}$ is a quasi-isomorphism.

Let us consider the case where
\(
	\mathcal{G} = \Right\Sheafhom{X}{\mathcal{F}(i), \omega_X[j-\dim X]}
\)
for some integers $i, j.$
In this case, we have
\[
	h \colon \mathcal{F} \longrightarrow \Right\Sheafhom{X}{\mathcal{F}(i),
	\omega_X[j-\dim X]}
\]
and its transpose morphism
\[
	\transpose h \colon \Right\Sheafhom{X}{\Right\Sheafhom{X}{\mathcal{F}(i),
	\omega_X[j-\dim X]}, \omega_X} \longrightarrow
	\Right\Sheafhom{X}{\mathcal{F},	\omega_X}.
\]
By abuse of notation, we denote the $(-i)$-th twist and the degree $(j-\dim X)$ shift of 
$\transpose h$ by the same symbol $\transpose h$:
\[
	\transpose h \colon \Right\Sheafhom{X}{\Right\Sheafhom{X}{\mathcal{F},
	\omega_X}, \omega_X} \longrightarrow
	\Right\Sheafhom{X}{\mathcal{F}(i),	\omega_X[j-\dim X]}
\]
By composing this morphism with
\[
	\canonical_{\mathcal{F}, \omega_X} \colon \mathcal{F} \longrightarrow
	\Right\Sheafhom{X}{\Right\Sheafhom{X}{\mathcal{F},
	\omega_X}, \omega_X},
\]
we obtain another morphism
\[
	\transpose h \circ \canonical_{\mathcal{F}, \omega_X} \colon 
	\mathcal{F} \longrightarrow \Right\Sheafhom{X}{\mathcal{F}(i), \omega_X[j-\dim X]}.
\]
\begin{dfn}\label{symmmor}
	Let $\mathcal{F}$ be a bounded complex of
	coherent $\strshf_X$-modules and $i, j$ some integers.
	A morphism 
	\[
		h \colon \mathcal{F} \longrightarrow \Right\Sheafhom{X}{\mathcal{F}(i),
		\omega_X[j-\dim X]}
	\]
	is said to be {\it symmetric} if it satisfies
	\[
		\transpose h  \circ \canonical_{\mathcal{F}, \omega_{X}}
		= h.
	\]
\end{dfn}

Let us introduce some notation on tuples of symmetric matrices. Let 
\[
	\Qdrcsp := k^{m+1} \otimes \Sym_2 k^{n+1}
\]
be the $k$-vector space of $(m+1)$-tuples of symmetric matrices of size $n+1$
with entries in $k.$ We write an ordered $k$-basis of the dual vector space $(k^{m+1})^\vee$
as $X_0, X_1, \dots, X_m.$ Then we can consider the basis as an ordered $k$-basis
of $\Cohomology{0}{\projsp^m, \strshf_{\projsp^m}(1)}.$
For an element $M=(M_0,M_1, \dots, M_m) \in \Qdrcsp,$
\[
	M(X) := X_0M_0 + X_1M_1 + \dots + X_mM_m
\]
is a symmetric matrix of size $n+1$ whose entries are $k$-linear forms in
$m+1$ variables $X_0 , X_1, \dots, X_m.$
We identify an $(m+1)$-tuple of symmetric matrices $M$ 
and a symmetric matrix $M(X)$ whose entries are $k$-linear forms in
$m+1$ variables. 
The {\it discriminant polynomial} of $M$ is defined by
\begin{align*}
	\disc (M) :=&  \det(M(X)) \\
	=&  \det(X_0M_0 + X_1M_1 + \dots + X_mM_m).
\end{align*}
If $\disc(M) \neq 0,$ the discriminant polynomial $\disc(M)$ is a homogeneous polynomial
of degree $n+1$ in $m+1$ variables. 
Let us define two subsets
\[
	\RedLoci \subset \NVLoci \subset \Qdrcsp = k^{m+1} \otimes \Sym_2 k^{n+1}.
\]

\begin{dfn}\label{nvloci}
Let $\NVLoci$ be the subset of $\Qdrcsp$ 
which consists of elements with non-zero discriminant polynomials.
The subset $\RedLoci \subset \NVLoci$ consists of elements
whose discriminant polynomials have no multiple factors over an algebraic closure of $k.$ 
\end{dfn}

\begin{dfn}\label{tctriple}
Let $\Triples$ be the set of equivalence classes of triples
 $(\mathcal{M}, \lambda, s),$ where
		\begin{itemize}
			\item $\mathcal{M}$ is a coherent $\strshf_{\projsp^m}$-module
				which is \ArithCM, pure of dimension $m-1,$ 
				and satisfies $\Cohomology{0}{\projsp^m, \mathcal{M}}=0$ and 
				$\dim \Cohomology{0}{\projsp^m, \mathcal{M}(1)} = n+1.$
			\item  $\lambda$ is a symmetric quasi-isomorphism
				\[
					\lambda \colon \mathcal{M} \isomarrow
					\Right\Sheafhom{\projsp^m}{\mathcal{M}(2-m), 
					\omega_{\projsp^m}[-m+1]}.
				\]
			\item $s= \{ s_0, s_1, \dots, s_n \}$  is an ordered $k$-basis of
				$\Cohomology{0}{\projsp^m, \mathcal{M}(1)}.$
			\end{itemize}
		Here, two triples $(\mathcal{M}, \lambda, s),
		(\mathcal{M}', \lambda', s')$ are said to be {\it equivalent} if
		there exists an isomorphism $\rho \colon \mathcal{M} \isomarrow \mathcal{M}'$
		of $\strshf_{\projsp^m}$-modules satisfying 
		\[
			\transpose \rho  \circ \lambda' \circ \rho
			= \lambda, \quad 
			\rho(s_i) = s'_i \quad (i = 0, 1, \dots, n).
		\]
\end{dfn}

Now we can state our first bijection in this paper.

\begin{thm}\label{bijection1}
	There exists a natural bijection 
	between $\NVLoci$ and $V_{m+1,n+1}.$
\end{thm}

The map from $\NVLoci$ to $V_{m+1,n+1}$ is denoted by
\[
	\phi_c \colon \NVLoci \longrightarrow V_{m+1, n+1},
\]
and the map from $V_{m+1, n+1}$ to $\NVLoci$ is denoted by
\[
	\psi_c \colon V_{m+1, n+1} \longrightarrow \NVLoci	.
\] 

The construction of the map $\psi_c$ is given in Subsection \ref{two-to-one},
and the construction of the map $\phi_c$ is given in Subsection \ref{one-to-two}.
Then we prove Theorem \ref{bijection1} in Subsection \ref{Proof}.

\begin{rmk}
	The maps $\phi_c, \psi_c$ depend on the choice of the quasi-isomorphism $c$
	fixed in the beginning of this section.
	(See Remark \ref{Change1} and Remark \ref{Change2} for details.)
\end{rmk}

\begin{rmk}
	We always assume $m \ge 2$ in this paper.
	This is because Proposition \ref{ACM-and-resolution} does not hold when $m=1.$
	It seems, however, that most of our arguments work well in the case of $m=1$ also.
	The case of $m=1$ and $n$ even is treated in \cite{Ishitsuka-Master}.
	For the study of the case of $m=1$ and its applications to the arithmetic of
	hyperelliptic curves and Fano schemes of lines, see \cite{Bhargava-Pencils}, \cite{Wang-Quadrics}.
\end{rmk}

\begin{rmk}
	The bijection of Theorem \ref{bijection1} says nothing about
	when there exists an $(m+1)$-tuple of symmetric matrices
	$M$ with $\disc(M) = f$ for a given homogeneous polynomial $f$ of degree $n+1$
	in $m+1$ variables. It would be an interesting problem in Arithmetic Invariant Theory 
	(cf. \cite[Theorem 23]{Bhargava-Pencils}).
\end{rmk}

\subsection{Construction of the map $\psi_c$}
\label{two-to-one}
In this subsection, we shall construct the map
\[
	\psi_c \colon V_{m+1, n+1} \longrightarrow \NVLoci.
\]

To avoid repeated arguments, we introduce the following notation and lemma.
Let $f \colon X \to Y$ be a proper morphism between schemes over $k$
of finite type. 
Let $\mathcal{F}$ (resp.\ $\mathcal{G}$) be a bounded complex of 
coherent $\strshf_X$-modules (resp.\ $\strshf_Y$-modules).
By Grothendieck duality, we have the following functorial quasi-isomorphism of 
complexes of coherent $\strshf_Y$-modules (\cite[Corollary 3.4 (c)]{Hartshorne-Residue}):
\[
	\GDual_f \colon \Right f_* \Right \Sheafhom{X}{\mathcal{F}, f^{!}\mathcal{G}}
	\isomarrow \Right \Sheafhom{Y}{\Right f_* \mathcal{F}, \mathcal{G}}.
\]

\begin{lem}\label{Cohomology}
	Let $\mathcal{M}$ be a coherent $\strshf_{\projsp^m}$-module 
	satisfying the following conditions:
	\begin{itemize}
		\item $\mathcal{M}$ is {\ArithCM} and pure of dimension $m-1,$ and
		\item there exists a quasi-isomorphism 
		\[
			\lambda \colon
			\mathcal{M} \isomarrow \Right\Sheafhom{\projsp^m}{\mathcal{M}(2-m),
		 \omega_{\projsp^m}[-m+1]}.
		 \]
	\end{itemize}
	Then we have
	\[
		\dim \Cohomology{i}{\projsp^m, \mathcal{M}(j)} = 
		\dim \Cohomology{m-i-1}{\projsp^m, \mathcal{M}(2-m-j)}
	\]
	for any $i, j.$
\end{lem}

\begin{prf}[Proof of Lemma \ref{Cohomology}]
		By Grothendieck duality for the structure morphism $f \colon \projsp^m \to \Spec k,$
	we obtain
	\begin{align*}
		\Right f_* (\mathcal{M}(j)) &\underset{\Right f_* \lambda}{\isomarrow}
		\Right f_* \Right \Sheafhom{\projsp^m}
		{\mathcal{M}(2-m-j), \omega_{\projsp^m}[-m+1]}\\
		&\underset{\GDual_f}{\isomarrow} \Right \Sheafhom{\Spec k}
		{\Right f_* \mathcal{M}(2-m-j), \strshf_{\Spec k}[-m+1]}.		
	\end{align*}
	Taking the cohomology, we have the desired equality. \qed
\end{prf}

Take an element $[(\mathcal{M}, \lambda, s)] \in \Triples.$ 
We have $\Cohomology{0}{\projsp^m, \mathcal{M}}=0$
by assumption. We also have $\Cohomology{m-1}{\projsp^m, \mathcal{M}(2-m)}=0$
by Lemma \ref{Cohomology}. Hence by Proposition \ref{ACM-and-resolution}, $\mathcal{M}$ admits 
a minimal graded locally free resolution of the following form for some $r \ge 0:$
\begin{equation}\label{Resolution1}
	\xymatrix{
		0 \ar[r] & {\displaystyle \bigoplus_{i=0}^r \strshf_{\projsp^m}(-2)} \ar[r]^M
		& {\displaystyle \bigoplus_{i=0}^r \strshf_{\projsp^m}(-1)} \ar[r]^(.67)p
		& \mathcal{M} \ar[r] & 0.
	}
\end{equation}

By the short exact sequence \Equref{Resolution1}, we obtain
\[
	\bigoplus_{i=0}^r \Cohomology{0}{\projsp^m, \strshf_{\projsp^m}} 
	\cong \Cohomology{0}{\projsp^m, \mathcal{M}(1)}.
\]
On the other hand, by our assumption on $\mathcal{M}$ 
(see Definition \ref{tctriple}), we have
\[
	\dim \Cohomology{0}{\projsp^m, \mathcal{M}(1)} = n+1.
\]
Hence we have $r = n.$ 
The map $M$ in \Equref{Resolution1} is identified with an $(m+1)$-tuple of
square matrices of size $n+1$ with entries in $k.$
However, these matrices are not necessarily symmetric.
We shall show that, by using $\lambda$ and $ s,$ 
we can choose a minimal graded locally free resolution of $\mathcal{M}$
so that we obtain a unique $(m+1)$-tuple of symmetric matrices.

First, we denote by $\{e_0, e_1, \dots, e_n\}$ the standard $\strshf_{\projsp^m}$-basis of
the middle term $\bigoplus_{i=0}^{n} \strshf_{\projsp^m}(-1)$ 
in \Equref{Resolution1}.
We put $\mathcal{G}:= \Kernel(p)$ and we denote the injection by
\[
	\iota \colon \mathcal{G} \hookrightarrow 
	\bigoplus_{i=0}^{n} \strshf_{\projsp^m}(-1)e_i.
 \]
By the short exact sequence \Equref{Resolution1},
$\mathcal{G}$ is isomorphic to $\bigoplus_{i=0}^{n} \strshf_{\projsp^m}(-2).$ 
We denote \Equref{Resolution1} by
\[
	\xymatrix{
		0 \ar[r] &
		\mathcal{G} \ar[r]^(.26){\iota} &
		{\displaystyle \bigoplus_{i=0}^{n} \strshf_{\projsp^m}(-1)e_i} \ar[r]^(.7)p &
		\mathcal{M} \ar[r] & 0.
	}
\]
We take a unique $p$ satisfying $p(e_i(1)) = s_i,$ 
where $\{e_0(1), e_1(1), \dots, e_n(1)\}$ denotes
the ordered $k$-basis of $\Cohomology{0}{\projsp^m, 
\bigoplus_{i=0}^{n} \strshf_{\projsp^m}e_i}$ corresponding
to $\{e_0, e_1, \dots, e_n\}.$

In order to simplify the notation, we define the functor 
\[
	\mathcal{F} \mapsto D\mathcal{F} := \Right\Sheafhom{\projsp^m}
	{\mathcal{F}(2-m), \omega_{\projsp^m}[-m+1]},
\]
and abbreviate the subindex of the canonical morphisms $\canonical.$
Then we have 
\[
	\lambda \colon \mathcal{M} \isomarrow D\mathcal{M}.
\]
In particular, $D\mathcal{M}$ is a sheaf.
It is easy to see from the definition that $D\mathcal{F}[-1]$ is a sheaf
if $\mathcal{F}$ is locally free.
We see that $Dh = \transpose h$ for any morphism $h \colon \mathcal{F} \longrightarrow
\mathcal{G}.$
Let us explain how to choose an appropriate isomorphism
\(
	\mathcal{G} \isomarrow \bigoplus_{i=0}^{n} \strshf_{\projsp^m}(-2).
\)

Applying the functor $D$ to the short exact sequence \Equref{Resolution1}, we obtain
\begin{equation}\label{Resolution2}
	\xymatrix{
		0 \ar[r] & D \left( {\displaystyle \bigoplus_{i=0}^{n}
		\strshf_{\projsp^m}(-1)e_i} \right)[-1] \ar[r]^(.7){\transpose \iota} &
		D{\mathcal{G}}[-1] \ar[r]^(.58){\delta} & 
		D\mathcal{M} \ar[r] & 0.
	}
\end{equation}

Now there exists a unique quasi-isomorphism
\[
	\rho \colon 
	\bigoplus_{i=0}^{n} \strshf_{\projsp^m}(-1)e_i \isomarrow D\mathcal{G}[-1]
\]
so that it satisfies
\[
	\delta \circ \rho = \lambda \circ p.
\]
This also gives a quasi-isomorphism
\[
	\transpose \rho \circ \canonical \colon 
	\mathcal{G} \isomarrow D(D\mathcal{G}) 
	\isomarrow D \left( \bigoplus_{i=0}^{n}	\strshf_{\projsp^m}(-1)e_i
	\right)[-1].
\]
Since an automorphism of $\mathcal{M}$ is uniquely lifted to
an automorphism of a pure minimal graded locally free resolution of
$\mathcal{M}$ by Lemma \ref{uniqueness},
we have the following commutative diagram:
\begin{equation}\label{Diagram1}
	\xymatrix{
		0 \ar[r] &
		\mathcal{G} \ar[r]^(.37){\iota} \ar[d]^{\xi}_{\isomvert} &
		{\displaystyle \bigoplus_{i=0}^{n} \strshf_{\projsp^m}(-1)e_i} 
		\ar[r]^(.7){p} \ar[d]^{\rho}_{\isomvert} &
		\mathcal{M} \ar[r] \ar[d]^{\lambda}_{\isomvert} & 0 \\
		0 \ar[r] &
		D \left( {\displaystyle \bigoplus_{i=0}^{n}	\strshf_{\projsp^m}(-1)e_i } \right)[-1]
		\ar[r]^(.65){\transpose \iota} &
		{D \mathcal{G}[-1]} 
		\ar[r]^(.55){\delta}&
		D\mathcal{M} \ar[r] & 0.
	}
\end{equation}
We shall show that the quasi-isomorphism $\xi$ in the diagram \Equref{Diagram1}
is equal to $\transpose \rho \circ \canonical.$
Applying the functor $D$ to the diagram \Equref{Diagram1}, 
we obtain the following commutative diagram:
\begin{equation}\label{Diagram2}
	\xymatrix{
		0 \ar[r] &
		D \left( {\displaystyle \bigoplus_{i=0}^{n}	\strshf_{\projsp^m}(-1)e_i } \right)[-1]
		\ar[r]^(.65){\transpose \iota} & D \mathcal{G}[-1] \ar[r]^(.55){\delta}  &
		D\mathcal{M} \ar[r] & 0 \\
		0 \ar[r] & \mathcal{G} \ar[r]^(.37){\iota} 
		\ar[u]_{\transpose \rho \circ \canonical}^{\isomvert}&
		{\displaystyle \bigoplus_{i=0}^{n} \strshf_{\projsp^m}(-1)e_i} 
		\ar[r]^(.7){p} \ar[u]_{\transpose \xi \circ \canonical}^{\isomvert}&
		\mathcal{M} \ar[r] \ar[u]_{\transpose \lambda 
		\circ \canonical}^{\isomvert}& 0.
	}
\end{equation}
Since the quasi-isomorphism $\lambda$ is symmetric by our assumption, we have
\(
	\transpose \lambda  \circ \canonical 
	= \lambda.
\)
We also have
\[
	\transpose \xi \circ \canonical =
	\rho \quad \Leftrightarrow \quad \xi = \transpose \rho \circ \canonical
\]
by the uniqueness of the lift of an automorphism of $\mathcal{M}$ 
(see Lemma \ref{uniqueness}).
Hence we can rewrite the diagram \Equref{Diagram1} as
\begin{equation}\label{Diagram3}
	\xymatrix{
		0 \ar[r] &
		D \left( {\displaystyle \bigoplus_{i=0}^{n}	\strshf_{\projsp^m}(-1)e_i } \right)[-1]
		 \ar[r]^(.6){\widetilde{M}} \ar@{=}[d] &
		{\displaystyle \bigoplus_{i=0}^{n} \strshf_{\projsp^m}(-1)e_i} 
		\ar[r]^(.65){p} \ar@{=}[d] &
		\mathcal{M} \ar[r] \ar[d]^{\lambda}_{\isomvert} & 0 \\
		0 \ar[r] &
		D \left( {\displaystyle \bigoplus_{i=0}^{n}	\strshf_{\projsp^m}(-1) e_i} \right)[-1]
		\ar[r]^(.6){ (*)} &
		{\displaystyle \bigoplus_{i=0}^{n} \strshf_{\projsp^m}(-1)e_i} 
		\ar[r]^(.65){q}&
		D\mathcal{M} \ar[r] & 0,
	}
\end{equation}
where we put $\widetilde{M} := \iota \circ \xi^{-1}$ and
$q := \delta \circ \rho = \lambda \circ p.$ 
The morphism $(*)$ in the diagram \Equref{Diagram3} is equal to 
$\rho^{-1} \circ \transpose \iota= \iota \circ \xi^{-1}.$ 
Since
\begin{align*}
	\canonical^{-1} \circ D \widetilde{M} &= 
	\canonical^{-1} \circ \transpose \xi^{-1} \circ \transpose \iota \\
	&= \rho^{-1} \circ \transpose \iota \\
	&= \iota \circ \xi^{-1} = \widetilde{M},
\end{align*}
we have $\canonical^{-1} \circ D \widetilde{M} = \widetilde{M}.$

Recall that we have fixed the quasi-isomorphism 
\[
	c \colon \omega_{\projsp^m} \isomarrow \strshf_{\projsp^m}(-m-1)[m]
\]
in the beginning of this section. By using $c,$ we obtain
\[
	D \left( \bigoplus_{i=0}^{n} \strshf_{\projsp^m}(-1)e_i \right)[-1]
	\isomarrow \bigoplus_{i=0}^{n} \strshf_{\projsp^m}(-2)
\]
and a symmetric matrix $M$ of size $n+1$
whose entries are in $\Cohomology{0}{\projsp^m, \strshf_{\projsp^m}(1)}.$ 
This matrix is identified with an $(m+1)$-tuple of symmetric matrices of size $n+1$
whose entries are in $k.$
This $(m+1)$-tuple is the desired one. 

Let us take another triple $(\mathcal{M}', \lambda', s')$ equivalent to 
$(\mathcal{M}, \lambda, s).$ There exists an isomorphism 
\[
	\rho \colon \mathcal{M} \isomarrow \mathcal{M}'
\]
satisfying $\rho(s_i) = s'_i$ and
$\transpose \rho \circ \lambda' \circ \rho = \lambda.$
Then we have the minimal graded locally free resolution of $\mathcal{M}'$
\[
	\xymatrix{
		0 \ar[r] &
		D \left( {\displaystyle \bigoplus_{i=0}^{n}	\strshf_{\projsp^m}(-1) e_i} \right)[-1]
		\ar[r]^(.6){\widetilde{M}'} &
		{\displaystyle \bigoplus_{i=0}^{n} \strshf_{\projsp^m}(-1)e_i} 
		\ar[r]^(.7){p'}&
		\mathcal{M}' \ar[r] & 0.
	}
\]
Since $\rho(s_i) = s'_i,$ we have $\rho \circ p = p'.$
Hence there exists an automorphism 
\[
	f \colon 
	D \left( {\displaystyle \bigoplus_{i=0}^{n} \strshf_{\projsp^m}(-1)e_i } \right)[-1]
	 \isomarrow 
	 D \left( {\displaystyle \bigoplus_{i=0}^{n}	\strshf_{\projsp^m}(-1)e_i } \right)[-1]
\] 
which makes the following diagram commute:
\[
	\xymatrix{
		0 \ar[r] &
		D \left( {\displaystyle \bigoplus_{i=0}^{n}	\strshf_{\projsp^m}(-1)e_i } \right)[-1]
		 \ar[r]^(.6){\widetilde{M}} \ar[d]_f^\vertisom &
		{\displaystyle \bigoplus_{i=0}^{n} \strshf_{\projsp^m}(-1)e_i} 
		\ar[r]^(.7){p} \ar@{=}[d]&
		\mathcal{M} \ar[r] \ar[d]^{\vertisom}_{\rho} & 0 \\
		0 \ar[r] &
		D \left( {\displaystyle \bigoplus_{i=0}^{n}	\strshf_{\projsp^m}(-1) e_i} \right)[-1]
		\ar[r]^(.6){\widetilde{M'}} &
		{\displaystyle \bigoplus_{i=0}^{n} \strshf_{\projsp^m}(-1)e_i} 
		\ar[r]^(.7){p'} &
		\mathcal{M}' \ar[r] & 0.
	}
\]
Dualizing this diagram, we have
\[
	\xymatrix{
		0 \ar[r] &
		D \left( {\displaystyle \bigoplus_{i=0}^{n}	\strshf_{\projsp^m}(-1)e_i } \right)[-1]
		 \ar[r]^(.6){\widetilde{M}} &
		{\displaystyle \bigoplus_{i=0}^{n} \strshf_{\projsp^m}(-1)e_i} 
		\ar[r]^(.65){q} &
		D\mathcal{M} \ar[r] & 0 \\
		0 \ar[r] &
		D \left( {\displaystyle \bigoplus_{i=0}^{n}	\strshf_{\projsp^m}(-1) e_i} \right)[-1]
		\ar[r]^(.6){\widetilde{M'}} \ar@{=}[u] &
		{\displaystyle \bigoplus_{i=0}^{n} \strshf_{\projsp^m}(-1)e_i} 
		\ar[r]^(.65){q'} \ar[u]_{\transpose f \circ \canonical}^\isomvert &
		D\mathcal{M}' \ar[r] \ar[u]^{\isomvert}_{\transpose \rho \circ \canonical} & 0.
	}
\]
Since $q' = \lambda' \circ p'$ and $\transpose \rho \circ \lambda' \circ \rho = \lambda,$
we have $\transpose \rho \circ q' = q.$
	On the other hand, we have
	\[
		\transpose \rho \circ q' = q \circ \transpose f. 
	\]
	Hence we have $q = q \circ \transpose f,$
	and by Lemma \ref{uniqueness}, we have $f=\identity.$
	Thus $M' = M,$ and any triple equivalent to $(\mathcal{M}, \lambda, s)$ 
	gives the same matrix $M.$

We put 
\[
	\psi_c \left([ (\mathcal{M}, \lambda, s)] \right) := M.
\] 
This finishes the construction of the map $\psi_c.$

\begin{rmk}\label{injection}
	By construction, the map $\psi_c$	is injective. 
	In fact, if 
	\[
		\psi_c \left( [(\mathcal{M}, \lambda, s)] \right) =
		\psi_c \left( [(\mathcal{M}', \lambda', s')] \right) = M,
	\]
	we have an isomorphism $\rho \colon \mathcal{M} \isomarrow \mathcal{M}'$ 
	and the following commutative diagram:
	\[
	\xymatrix{
		0 \ar[r] &
		D \left( {\displaystyle \bigoplus_{i=0}^{n}	\strshf_{\projsp^m}(-1)e_i } \right)[-1]
		 \ar[r]^(.6){\widetilde{M}} \ar@{=}[d] &
		{\displaystyle \bigoplus_{i=0}^{n} \strshf_{\projsp^m}(-1)e_i} 
		\ar[r]^(.7){p} \ar@{=}[d]&
		\mathcal{M} \ar[r] \ar[d]^{\vertisom}_{\rho} & 0 \\
		0 \ar[r] &
		D \left( {\displaystyle \bigoplus_{i=0}^{n}	\strshf_{\projsp^m}(-1) e_i} \right)[-1]
		\ar[r]^(.6){\widetilde{M}} &
		{\displaystyle \bigoplus_{i=0}^{n} \strshf_{\projsp^m}(-1)e_i} 
		\ar[r]^(.7){p'} &
		\mathcal{M}' \ar[r] & 0.
	}
	\]
	We see that the isomorphism $\rho$ satisfies
	\[
		\transpose \rho  \circ \lambda' \circ \rho
		= \lambda
	\]
	because $\canonical^{-1} \circ D\widetilde{M} = \widetilde{M}.$ 
	Hence two triples $(\mathcal{M}, \lambda, s), (\mathcal{M}', \lambda', s')$ are 
	equivalent to each other.
\end{rmk}

\begin{rmk}\label{Change1}
	Note that we use the quasi-isomorphism $c$ only in the last step.
	If one uses $ac$ for some $a \in k^\times$ instead of $c,$
	we obtain another map $\psi_{ac}.$
	(Note that, since $\omega_{\projsp^m}$ and 
	$\strshf_{\projsp^m}(-m-1)$ are represented by 
	some degree shift of a line bundle on $\projsp^m,$
	any other quasi-isomorphism 
	\[
		c' \colon \omega_{\projsp^m} \isomarrow
		\strshf_{\projsp^m}(-m-1)[-m]
	\]
	can be written as $c' = ac$ for some $a \in k^\times.$)
	These two maps $\psi_c$ and $\psi_{ac}$ satisfy
	\[
		\psi_{ac}\left([(\mathcal{M}, \lambda, s)]\right) = 
		a^{-1}\psi_c \left([(\mathcal{M}, \lambda, s)]\right).
	\] 
	This can be seen easily from the diagram \Equref{Diagram3}.
\end{rmk}

\subsection{Construction of the map $\phi_c$}\label{one-to-two}
In this section, we construct the map
\[
	\phi_c \colon \NVLoci \longrightarrow \Triples.
\]
We go backward in the steps of the construction of $\psi_c$ in the previous subsection.
Let $M \in \NVLoci$ be a symmetric matrix of size $n+1$
such that the entries of $M$ are $k$-linear forms in $m+1$ variables $X_0, X_1, \dots, X_m$
and $M$ satisfies $\disc(M) \neq 0.$
Then, we immediately obtain the following short exact sequence
\begin{equation}\label{Resolution3}
	\xymatrix{
		0 \ar[r] & {\displaystyle \bigoplus_{i=0}^{n} \strshf_{\projsp^m}(-2)} \ar[r]^{M}
		& {\displaystyle \bigoplus_{i=0}^{n} \strshf_{\projsp^m}(-1)} \ar[r]^(.7){p}
		& \mathcal{M} \ar[r] & 0,
	}
\end{equation}
where we put $\mathcal{M} := \Cokernel(M).$
By Proposition \ref{ACM-and-resolution},
the cokernel $\mathcal{M}$ of $M$ is \ArithCM, pure of dimension $m-1$ and satisfies
\[
	\Cohomology{0}{\projsp^m, \mathcal{M}} =
	\Cohomology{m-1}{\projsp^m, \mathcal{M}(2-m)} = 0.
\]
So $\mathcal{M}$ is $1$-regular and $\mathcal{M}(1)$ is generated by its global sections.
By the short exact sequence \Equref{Resolution3}, the morphism
\[
	\Cohomology{0}{\projsp^m, \bigoplus_{i=0}^n \strshf_{\projsp^m}e_i} \longrightarrow
	\Cohomology{0}{\projsp^m, \mathcal{M}(1)}
\]
is an isomorphism. In particular, we have $\dim \Cohomology{0}{\projsp^m, 
\mathcal{M}(1)}=n+1.$
Take $s = \{s_0, s_1, \dots, s_n\}$ to be the image of the standard basis
\(
	\left\{ e_o, e_1, \dots, e_n \right\}
\)
of $\bigoplus_{i=0}^n \strshf_{\projsp^m}e_i.$
Since $M$ is a symmetric matrix, we have the same resolution for 
\(
	\Sheafext{1}{\projsp^m}{\mathcal{M}(2-m), \omega_{\projsp^m}[-m+1]}
\)
and 
\[
	\Sheafext{i}{\projsp^m}{\mathcal{M}(2-m), \omega_{\projsp^m}[-m+1]} = 0
\]
for $i \neq 1.$ Hence we obtain a quasi-isomorphism
\[
	\lambda \colon \mathcal{M} \isomarrow \Right\Sheafhom{\projsp^m}
	{\mathcal{M}(2-m), \omega_{\projsp^m}[-m+1]}.
\] 

We shall show the quasi-isomorphism $\lambda$ obtained above is symmetric.
(See Definition \ref{symmmor}.) 
By using $c,$ we modify the exact sequence \Equref{Resolution3} to
\begin{equation}\label{Resolution4}
	\xymatrix{
		0 \ar[r] & D \left( {\displaystyle \bigoplus_{i=0}^{n} 
		\strshf_{\projsp^m}(-1)e_i} \right)[-1] \ar[r]^(.6){\widetilde{M}}
		& {\displaystyle \bigoplus_{i=0}^{n} \strshf_{\projsp^m}(-1)e_i} \ar[r]^(.7){p}
		& \mathcal{M} \ar[r] & 0.
	}
\end{equation}
So we obtain two commutative diagrams
\begin{equation}\label{Diagram4}
	\xymatrix{
		0 \ar[r] &
		D \left( {\displaystyle \bigoplus_{i=0}^{n}	\strshf_{\projsp^m}(-1)e_i } \right)[-1]
		 \ar[r]^(.6){\widetilde{M}} \ar@{=}[d] &
		{\displaystyle \bigoplus_{i=0}^{n} \strshf_{\projsp^m}(-1)e_i} 
		\ar[r]^(.68){p} \ar@{=}[d] &
		\mathcal{M} \ar[r] \ar[d]^{\lambda}_{\isomvert} & 0 \\
		0 \ar[r] &
		D \left( {\displaystyle \bigoplus_{i=0}^{n}	\strshf_{\projsp^m}(-1)e_i } \right)[-1]
		\ar[r]^(.6){\transpose \widetilde{M} = \widetilde{M}} &
		{\displaystyle \bigoplus_{i=0}^{n} \strshf_{\projsp^m}(-1)e_i} 
		\ar[r]^(.68){\delta}&
		D\mathcal{M} \ar[r] & 0
	}
\end{equation}
and 
\[
	\xymatrix{
		0 \ar[r] &
		D \left( {\displaystyle \bigoplus_{i=0}^{n}	\strshf_{\projsp^m}(-1)e_i } \right)[-1]
		 \ar[r]^(.6){\widetilde{M}} \ar@{=}[d] &
		{\displaystyle \bigoplus_{i=0}^{n} \strshf_{\projsp^m}(-1)e_i} 
		\ar[r]^(.68){p} \ar@{=}[d] &
		\mathcal{M} \ar[r]  
		\ar[d]^{\transpose \lambda  \circ 
		\canonical}_{\isomvert}& 0 \\
		0 \ar[r] &
		D \left( {\displaystyle \bigoplus_{i=0}^{n}	\strshf_{\projsp^m}(-1)e_i } \right)[-1]
		\ar[r]^(.6){\transpose \widetilde{M} = \widetilde{M}} &
		{\displaystyle \bigoplus_{i=0}^{n} \strshf_{\projsp^m}(-1)e_i} 
		\ar[r]^(.68){\delta}&
		D\mathcal{M} \ar[r] 
		& 0.
	}
\]
Comparing these diagrams, we conclude that $\lambda$ is symmetric.

\begin{rmk}\label{Change2}
	The map $\phi_c$ depends on the choice of the quasi-isomorphism $c.$
	If we use $c' = ac$ instead of $c$ for $a \in k^\times,$ 
	we obtain $\phi_{ac}(M) = [(\mathcal{M}, a\lambda, s)]$ instead of 
	$\phi_c(M) = [(\mathcal{M}, \lambda, s)].$
	In fact, the coherent $\strshf_{\projsp^m}$-module 
	$\mathcal{M}$ and the morphism $p$ does not change,
	so the ordered basis $s$ is not affected. On the other hand, 
	the morphism $\delta$ changes to
	$a\delta.$ Because $\lambda$ satisfies 
	$\lambda \circ p = \delta,$ we see that
	$\lambda$ changes to $a\lambda.$
	We also see that $\phi_c(aM) = [(\mathcal{M}, a^{-1}\lambda,
	s)]$ by a similar argument.
	(Note that the triple $(\mathcal{M}, \lambda, bs)$ is equivalent to
	the triple $(\mathcal{M}, b^2\lambda, s)$ for $b \in k^\times.$)
\end{rmk}

\subsection{End of the proof of Theorem \ref{bijection1}}\label{Proof}
We shall show that the maps $\phi_c, \psi_c$ constructed 
in Subsection \ref{two-to-one} and Subsection \ref{one-to-two} are 
inverses to each other.
By Remark \ref{injection}, the map $\psi_c$ is injective. 
Hence we only have to show that 
\[
	\psi_c \circ \phi_c \colon \NVLoci \longrightarrow \NVLoci
\]
is the identity map.

Take a symmetric matrix $M \in \NVLoci.$
Then we obtain a triple $(\mathcal{M}, \lambda, s)$
and the diagram \Equref{Diagram4} by using the fixed quasi-isomorphism $c.$
On the other hand, from the triple, we obtain the matrix $M'$
and the resolution
\[
	\xymatrix{
		0 \ar[r] & 
		D \left( {\displaystyle \bigoplus_{i=0}^{n}	\strshf_{\projsp^m}(-1) e_i} \right)[-1]
		 \ar[r]^(.6){\widetilde{M'}} &
		{\displaystyle \bigoplus_{i=0}^{n} \strshf_{\projsp^m}(-1) e_i} 
		\ar[r]^(.7){p'} &
		\mathcal{M} \ar[r] & 0.
	}
\]
In other words, we have $\psi_c([(\mathcal{M}, \lambda, s)]) = M'.$
We must check that $M=M'.$

Since both of the morphisms $p$ and $p'$ are determined 
by the same ordered basis $s = \{s_0, s_1, \dots, s_n\}$ of 
$\Cohomology{0}{\projsp^m, \mathcal{M}(1)},$
they are the same morphism. 
There exists an automorphism 
\[
	f \colon D \left( {\displaystyle \bigoplus_{i=0}^{n} 
	\strshf_{\projsp^m}(-1)e_i } \right)[-1]
	\isomarrow D \left( {\displaystyle \bigoplus_{i=0}^{n} 
	\strshf_{\projsp^m}(-1)e_i } \right)[-1]
\]
such that the following diagram is commutative:
\[
	\xymatrix{
		0 \ar[r] &
		D \left( {\displaystyle \bigoplus_{i=0}^{n}	\strshf_{\projsp^m}(-1)e_i } \right)[-1]
		 \ar[r]^(.6){\widetilde{M}} \ar[d]_f^{\vertisom} &
		{\displaystyle \bigoplus_{i=0}^{n} \strshf_{\projsp^m}(-1)e_i} 
		\ar[r]^(.7){p} \ar@{=}[d] &
		\mathcal{M} \ar[r]  
		\ar@{=}[d] & 0 \\
		0 \ar[r] &
		D \left( {\displaystyle \bigoplus_{i=0}^{n}	\strshf_{\projsp^m}(-1)e_i } \right)[-1]
		\ar[r]^(.6){\widetilde{M'}} &
		{\displaystyle \bigoplus_{i=0}^{n} \strshf_{\projsp^m}(-1)e_i} 
		\ar[r]^(.7){p' = p}&
		\mathcal{M} \ar[r] 
		& 0.
	}
\]
Dualizing this diagram, we obtain
\[
	\xymatrix{
		0 \ar[r] &
		D \left( {\displaystyle \bigoplus_{i=0}^{n}	\strshf_{\projsp^m}(-1)e_i } \right)[-1]
		 \ar[r]^(.6){\widetilde{M}} \ar@{=}[d] &
		{\displaystyle \bigoplus_{i=0}^{n} \strshf_{\projsp^m}(-1)e_i} 
		\ar[r]^(.68){q}  &
		D\mathcal{M} \ar[r]  
		\ar@{=}[d] & 0 \\
		0 \ar[r] &
		D \left( {\displaystyle \bigoplus_{i=0}^{n}	\strshf_{\projsp^m}(-1)e_i } \right)[-1]
		\ar[r]^(.6){\widetilde{M'}} &
		{\displaystyle \bigoplus_{i=0}^{n} \strshf_{\projsp^m}(-1)e_i} 
		\ar[r]^(.68){q} \ar[u]_{\transpose f  \circ
		\canonical}^{\isomvert}&
		D\mathcal{M} \ar[r] 
		& 0.
	}
\]
Since $D\mathcal{M}(1) \cong \mathcal{M}(1)$ is
generated by its global sections, we have $f = \identity.$
This completes the proof of Theorem \ref{bijection1}.

\section{Actions of $\GenLin_{m+1}(k) \times \GenLin_{n+1}(k)$ 
on $\NVLoci$ and $\Triples$ }\label{GenLin}
In this section, we introduce actions of $\GenLin_{m+1}(k) \times \GenLin_{n+1}(k)$ 
on $\NVLoci$ and  $\Triples$
so that the maps $\phi_c, \psi_c$ defined in Section \ref{bijection} are equivariant.
We use the same notation as in Section \ref{bijection}.

Take $M=(M_0,M_1, \dots, M_m) \in \Qdrcsp.$ 
For $A = (a_{i,j})_{0 \le i,j \le m} \in 
\GenLin_{m+1}(k)$ and $P \in \GenLin_{n+1}(k),$ we define
\[
	M \cdot (A,P) = 
	\left( \sum_{i=0}^m a_{i,0} \transpose P M_i P, \;
	 \sum_{i=0}^m a_{i,1} \transpose P M_i P, \; \dots, \;
	 \sum_{i=0}^m a_{i,m} \transpose P M_i P\right)
\]
and
\[
	M \cdot A := M \cdot (A, I_{n+1}),
\]
where $I_r$ is the identity matrix of size $r.$
We also write $\transpose P M P := M \cdot (I_{m+1}, P).$
For $X = (X_0, X_1, \dots, X_m) \in (k^{m+1})^\vee := 
\bigoplus_{i=0}^m \Homomorphism{k}{k, k},$
we define
\begin{align*}
	AX &:=   \left( 
	\sum_{j=0}^m a_{0,j}X_j, 
	\sum_{j=0}^m a_{1,j}X_j, \dots
	\sum_{j=0}^m a_{m,j}X_j \right).
\end{align*}
Then we have
\[
	\left( M \cdot (A, P) \right)(X) = \transpose P M(AX) P. 
\]
To study the corresponding action of $\GenLin_{m+1}(k) \times \GenLin_{n+1}(k)$
on the set $\Triples,$ it is enough to check the actions of $\GenLin_{m+1}(k)$ 
and $\GenLin_{n+1}(k)$ separately.

First we check the action of $\GenLin_{n+1}(k).$
Put $\phi_c(\transpose PMP) = [(\mathcal{M}', \lambda', s')] \in \Triples.$ 
We shall show $(\mathcal{M}', \lambda', s')$ is equivalent to 
$(\mathcal{M}, \lambda, s \transpose P^{-1}).$
Here, for an ordered $k$-basis $s=\{s_0, s_1, \dots, s_n\}$ 
and $P = (b_{i,j}) \in \GenLin_{n+1}(k),$
we define
\[
	s P := \left\{ \sum_{i = 0}^n s_ib_{i, 0}, \;
	\sum_{i=0}^n s_i b_{i, 1}, \; \dots, \; 
	\sum_{i=0}^n s_i b_{i, n} \right\}.
\]
The sequence \Equref{Resolution4} attached to $\transpose P M P$ is
\[
	\xymatrix{
		0 \ar[r] &
		D \left( {\displaystyle \bigoplus_{i=0}^{n}	\strshf_{\projsp^m}(-1) e_i} \right)[-1]
		\ar[r]^(.6){\transpose P\widetilde{M}P} &
		{\displaystyle \bigoplus_{i=0}^{n} \strshf_{\projsp^m}(-1)e_i} 
		\ar[r]^(.7){p'} &
		\mathcal{M}' \ar[r] & 0.
	}
\]
We combine the sequences \Equref{Resolution4} for
$M$ and $\transpose P M P$ as follows:
\[
	\xymatrix{
		0 \ar[r] &
		D \left( {\displaystyle \bigoplus_{i=0}^{n}	\strshf_{\projsp^m}(-1)e_i } \right)[-1]
		 \ar[r]^(.6){\widetilde{M}} \ar[d]^{P^{-1}} &
		{\displaystyle \bigoplus_{i=0}^{n} \strshf_{\projsp^m}(-1)e_i} 
		\ar[r]^(.7){p} \ar[d]^{\transpose P}  &
		\mathcal{M} \ar[r] \ar[d]^{\vertisom}_{\rho} & 0 \\
		0 \ar[r] &
		D \left( {\displaystyle \bigoplus_{i=0}^{n}	\strshf_{\projsp^m}(-1) e_i} \right)[-1]
		\ar[r]^(.6){\transpose P\widetilde{M}P} &
		{\displaystyle \bigoplus_{i=0}^{n} \strshf_{\projsp^m}(-1)e_i} 
		\ar[r]^(.7){p'} &
		\mathcal{M}' \ar[r] & 0.
	}
\]
Here an isomorphism $\rho$ exists because the left square commutes.
The isomorphism $\rho$ satisfies 
\[
	\rho \circ p(e_i) = p'(e_i \transpose P^{-1}).
\]
This tells us that the ordered $k$-basis $s'$ corresponding to
$\transpose PMP$ is equal to $s \transpose P^{-1}.$

Moreover, applying $D$ to this diagram, we find
\[
	\xymatrix{
		0 \ar[r] &
		D \left( {\displaystyle \bigoplus_{i=0}^{n}	\strshf_{\projsp^m}(-1)e_i } \right)[-1]
		 \ar[r]^(.6){\widetilde{M}} \ar[d]^{P^{-1}} &
		{\displaystyle \bigoplus_{i=0}^{n} \strshf_{\projsp^m}(-1)e_i} 
		\ar[r]^(.7){\delta} \ar[d]^{\transpose P}  &
		D\mathcal{M} \ar[r]  
		\ar[d]_{\transpose \rho^{-1}  \circ \canonical}^{\vertisom} & 0 \\
		0 \ar[r] &
		D \left( {\displaystyle \bigoplus_{i=0}^{n}	\strshf_{\projsp^m}(-1)e_i } \right)[-1]
		\ar[r]^(.6){\transpose P\widetilde{M}P} &
		{\displaystyle \bigoplus_{i=0}^{n} \strshf_{\projsp^m}(-1)e_i} 
		\ar[r]^(.7){\delta'} &
		D\mathcal{M} \ar[r] & 0.
	}
\]
With these diagrams, we have the following equalities:
\begin{alignat*}{2}
	\rho \circ p &= p' \circ \transpose P, &\qquad 
	\lambda \circ p &= \delta, \\
	\transpose \rho^{-1} \circ \delta &=
	 \delta' \circ \transpose P, & \qquad 
	\lambda' \circ p' &= \delta', 
\end{alignat*}
and since the left squares in two diagrams are the same, we see that
\[
	\transpose \rho  \circ \lambda' \circ \rho
	= \lambda.
\]

Next we examine the action of $\GenLin_{m+1}(k).$
Put $\phi_c (M \cdot A)=[(\mathcal{M}'', \lambda'', s'')].$
We write 
\[
	\nu_A \colon \projsp^m \isomarrow \projsp^m \quad ; \quad 
	[u_0 : u_1 : \dots : u_m] \mapsto 
	\left[ \sum_{j=0}^{m} a_{0,j}u_j : \sum_{j=0}^{m} a_{1,j}u_j : \dots : 
	\sum_{j=0}^{m} a_{m,j}u_j \right]
\]
the projective automorphism induced by $A = (a_{i,j})_{0 \le i, j \le m}.$
The triple $(\mathcal{M}'', \lambda'', s'')$ is equivalent to 
$(\nu_A^*\mathcal{M}, \nu_A^*\lambda, \nu_A^*s).$
In fact, if we apply $\nu_A^*$ to the sequence \Equref{Resolution3},
we have
\[
	\xymatrix{
		0 \ar[r] & D \left( {\displaystyle \bigoplus_{i=0}^{n} 
		\strshf_{\projsp^m}(-1) e_i} \right)[-1]
		\ar[r]^(.6){\widetilde{\nu_A^*M}} &
		{\displaystyle \bigoplus_{i=0}^{n} \strshf_{\projsp^m}(-1)e_i} 
		\ar[r]^(.65){\nu_A^* p} &
		\nu_A^* \mathcal{M} \ar[r] & 0.
	}
\]
Since 
\[
	(M \cdot A)(X) = M(X)\cdot A =  M (AX) = M(\nu_A^*X) =(\nu_A^* M)(X),
\] 
we have $\widetilde{\nu_A^*M} = \widetilde{M \cdot A} = \widetilde{M} \cdot A.$ 
Hence there is an isomorphism 
\[
	\rho \colon \mathcal{M}'' \isomarrow \nu_A^* \mathcal{M}
\]
which makes the diagram
\[
	\xymatrix{
		0 \ar[r] &
		D\left( {\displaystyle \bigoplus_{i=0}^{n}	\strshf_{\projsp^m}(-1) e_i} \right)[-1] 
		 \ar[r]^(.6){\widetilde{M} \cdot A} \ar@{=}[d] &
		{\displaystyle \bigoplus_{i=0}^{n} \strshf_{\projsp^m}(-1)e_i} 
		\ar[r]^(.7){p''} \ar@{=}[d]  &
		\mathcal{M}'' \ar[r] \ar[d]^{\vertisom}_{\rho} & 0 \\
		0 \ar[r] &
		D \left( {\displaystyle \bigoplus_{i=0}^{n} \strshf_{\projsp^m}(-1) e_i} \right)[-1]
		\ar[r]^(.6){\widetilde{M} \cdot A} &
		{\displaystyle \bigoplus_{i=0}^{n} \strshf_{\projsp^m}(-1)e_i} 
		\ar[r]^(.7){\nu_A^* p} &
		\nu_A^* \mathcal{M} \ar[r] & 0
	}
\]
commute.
In the similar way to the case of $\GenLin_{n+1}(k)$-action, we find 
\[
	\transpose \rho  \circ \nu_A^*\lambda \circ \rho =  \lambda''.
\]
This shows the required equivalence.
The actions of $\GenLin_{m+1}(k)$ and $ \GenLin_{n+1}(k)$ commute,
so $\nu_A^* (s \transpose P^{-1}) = (\nu_A^*s) \transpose P^{-1}.$
We write it as $\nu_A^*s \transpose P^{-1}.$
With these arguments, we find
\[
	\phi_c(M \cdot (A, P)) = [(\nu_A^* \mathcal{M}, \nu_A^* \lambda,
	\nu_A^*s  \transpose P^{-1} )].
\]
Hence we have defined the action of $\GenLin_{m+1}(k) \times \GenLin_{n+1}(k)$ on
$\Triples$ as 
\[
	[(\mathcal{M}, \lambda, s)] \cdot (A, P) :=
	[(\nu_A^* \mathcal{M}, \nu_A^* \lambda, \nu_A^*s  \transpose P^{-1} )].
\]
Then we conclude
\[
	\phi_{c}(M \cdot (A, P)) = \phi_c(M) \cdot (A, P).
\]

\begin{rmk}
	The group $\GenLin_{m+1}(k)$ naturally acts on the space $\Homomorphism{\projsp^m}
	{\omega_{\projsp^m}, \strshf_{\projsp^m}(-m-1)[m]}$ via 
	the inverse of the determinant, i.e.
	\[
		\nu_A^*c = \det(A)^{-1} c.
	\]
	We must take this action in consideration if we would like to 
	make the bijections free from the choice of the quasi-isomorphism
	\[
		c \colon \omega_{\projsp^m} \isomarrow \strshf_{\projsp^m}(-m-1)[m].
	\]
	However, if we modify the action of $\GenLin_{m+1}(k) \times \GenLin_{n+1}(k)$ 
	on $\NVLoci$ as
	\[
		M * (A, P) := \det (A)^{-1} M \cdot (A, P),
	\]
	we obtain
	\begin{align*}
		\phi_{c \cdot A}(M * (A, P)) &= \phi_{\det(A)^{-1} c} 
		(\det(A)^{-1}M \cdot (A, P)) \\
		&= \phi_c(M \cdot (A, P)) \\
		&= \phi_c (M) \cdot (A, P)
	\end{align*}
	from Remark \ref{Change2}.
	Hence we can justify our arguments.
\end{rmk}

As a special case, let us consider $(A, P) = (aI_{m+1}, I_{n+1})$ for some $a \in k^\times.$ 
We put $\phi_c(M) = [(\mathcal{M}, a^{-1}\lambda, s)].$ 
Then by Remark \ref{Change2}, we have
\[
	\phi_c(M \cdot (A, P)) = \phi_c(aM) = [(\mathcal{M}, a^{-1}\lambda, s)].
\]
To consider bijections between orbits, we introduce the following 
equivalence relations $\sim_1$ and $\sim_2.$
\begin{dfn}\label{sims}
	Let $(\mathcal{M}, \lambda)$ be a pair satisfying the following conditions:
	\begin{itemize}
		\item $\mathcal{M}$ is a coherent $\strshf_{\projsp^m}$-module
		which is {\ArithCM}, pure of dimension $m-1,$ and satisfies
		$\Cohomology{0}{\projsp^m,\mathcal{M}} = 0$ and
		$\dim_k \Cohomology{0}{\projsp^m, \mathcal{M}(1)} = n+1.$
		\item $\lambda$ is a symmetric quasi-isomorphism
		\[
			\lambda \colon \mathcal{M} \isomarrow
			\Right\Sheafhom{\projsp^m}{\mathcal{M}(2-m), \omega_{\projsp^m}[-m+1]}.
		\] 
	\end{itemize}
	Let $(\mathcal{M}', \lambda')$ be another pair satisfying the same
	conditions  as $(\mathcal{M}, \lambda).$
	\begin{enumerate}
	\item We write $(\mathcal{M}, \lambda) \sim_1 (\mathcal{M}', \lambda')$
	if there exists an automorphism 
	$\rho \colon \mathcal{M} \isomarrow \mathcal{M}'$ satisfying 
	\[
		\transpose \rho \circ \lambda' \circ \rho = u \lambda 
	\]
	for some $u \in k^\times.$
	
	\item We write $(\mathcal{M}, \lambda) \sim_2 (\mathcal{M}', \lambda')$ 
	if there exist a matrix $A \in \GenLin_{m+1}(k)$ and an isomorphism
	$\rho \colon \nu_A^* \mathcal{M} \isomarrow \mathcal{M}'$ satisfying 
	\[
		\transpose \rho \circ \lambda' \circ \rho =\nu_A^* \lambda.
	\] 
	\end{enumerate}
\end{dfn}
\begin{cor}\label{Bij-cor1}
	There exists a natural bijection between the following two sets.
	\begin{itemize}
		\item\label{matrices1-1}
		The set of $(k^\times I_{m+1}) \times \GenLin_{n+1}(k)$-orbits of 
		$(m+1)$-tuples of symmetric matrices 
		$M = (M_0,M_1,\dots,M_m)$ of size $n+1$ with entries in $k$
		satisfying $\disc (M) \neq 0.$  
		\item\label{sheaves1-1}
		The set of equivalence classes of pairs $(\mathcal{M}, \lambda)$ 
		satisfying the conditions of Definition \ref{sims}
		with respect to the equivalence relation $\sim_1.$
	\end{itemize}
\end{cor}

Meanwhile, the orbits of the whole group $\GenLin_{m+1}(k) \times \GenLin_{n+1}(k)$ 
gives another bijection:
\begin{cor}\label{Bij-cor2}
	There exists a natural bijection between the following two sets.
	\begin{itemize}
		\item\label{matrices1-2}
		The set of $\GenLin_{m+1}(k) \times \GenLin_{n+1}(k)$-orbits of 
		$(m+1)$-tuples of symmetric matrices 
		$M = (M_0,M_1,\dots,M_m)$ of size $n+1$ with entries in $k$
		satisfying $\disc (M) \neq 0.$  
		\item\label{sheaves1-2}
		The set of equivalence classes of pairs $(\mathcal{M}, \lambda)$
		satisfying the conditions of Definition \ref{sims}
		with respect to the equivalence relation $\sim_2.$
	\end{itemize}
\end{cor}

\section{Description of endomorphisms and fibers}\label{fibers}
We follow the notation used in Section \ref{bijection} and Section \ref{GenLin}.
In this section, we study the set of matrices $M \in \NVLoci$
such that the pair $(\mathcal{M}, \lambda)$ 
defined by triples $\phi_c(M) = [(\mathcal{M}, \lambda, s)]$ are
equivalent to a fixed pair $(\mathcal{M}', \lambda')$
with respect to the equivalence relation $\sim_1$ (resp.\ $\sim_2$).
To describe in algebraic terms, we introduce the following 
two equivalence relations on symmetric quasi-isomorphisms.

Let $M \in \NVLoci$ be an element.
We put $\phi_c(M) = [(\mathcal{M}, \lambda, s)]$ and
write the closed subscheme defined by $\disc(M)$ by
\[
	\Delta_M:= \Supp(\mathcal{M}) 
	= (\disc(M)=0) \subset \projsp^m.
\]
We start with discussing symmetric homomorphisms from $\mathcal{M}$ to 
\[
	D\mathcal{M} := \Right\Sheafhom{\projsp^m}{\mathcal{M}(2-m), 
	\omega_{\projsp^m}[-m+1]}.
\]

We put $L_0 := \Endomorphism{\projsp^m}{\mathcal{M}}.$
Then $L_0$ is a finite-dimensional $k$-algebra,
and the $k$-vector space 
\(
	\Homomorphism{\projsp^m}{\mathcal{M}, 
	D\mathcal{M}} 
\)
has two ways of simply transitive actions of $L_0.$ Namely, for $l \in L_0,$ we have an action
\[
	\eta \mapsto \eta \circ l
\]
and another action
\[
	\eta \mapsto \transpose l  \circ \eta.
\]
Assume that $\eta$ is symmetric. Then we have
\(
	\transpose (\eta \circ l) = \transpose l \circ \transpose \eta = \transpose l  \circ \eta.
\)
Henceforth $\eta \circ l$ is symmetric if and only if
\begin{equation}\label{Symmetricity}
	\transpose l  \circ \eta
	= \eta \circ l.
\end{equation}

Recall that $\lambda$ is a symmetric quasi-isomorphism.
Let $L$ denote the subspace of $L_0$ consisting of elements $l \in L_0$
satisfying 
\[
	\transpose l \circ \lambda = \lambda \circ l.
\]
Then the map
\[
	L \longrightarrow \Homomorphism{\projsp^m}{\mathcal{M}, 
	D\mathcal{M}} \quad ; \quad
	l \mapsto \lambda \circ l
\]
gives a bijection between the subspace $L \subset L_0$ and
the subspace of symmetric homomorphisms.
In terms of the matrix algebra $\Matrix_{n+1}(k),$ we can describe $L_0$ and $L$ as follows.

\begin{prop}\label{End1}
	The $k$-algebra $L_0$ is isomorphic to
	\begin{align}\label{endom1}
		\left\{ (P, P') \in \Matrix_{n+1}(k) \times \Matrix_{n+1}(k) \; \relmid \;
		\transpose P M = M P' \right\}
	\end{align}
	where the product of $\Matrix_{n+1}(k) \times \Matrix_{n+1}(k)$ is defined by
	\[
		(P_1, P'_1) \circ (P_2, P'_2) = (P_2 P_1, P'_1 P'_2).
	\]
	The anti-homomorphism $(P, P') \mapsto (P', P)$ gives 
	an anti-endomorphism $\sigma$ of $L_0$ as a $k$-algebra, 
	and the fixed part $L^\sigma_0$ is identified with $L.$
\end{prop}

\begin{rmk}
	The fixed part $L = L^\sigma_0$ is not a $k$-subalgebra of $L_0$ in general.
	In other words, for two endomorphisms $l$ and $l'$ such that $\lambda \circ l$ and
	$\lambda \circ l'$ are symmetric,
	it is not always the case that $\lambda \circ l \circ l'$ is symmetric.	
\end{rmk}

\begin{prf}
	We first show the existence of an inclusion 
	$L_0 \hookrightarrow \Matrix_{n+1}(k) \times \Matrix_{n+1}(k).$
	Take $l \in L_0.$  Then, since an endomorphism of $\mathcal{M}$ is
	uniquely lifted to an endomorphism of a pure minimal graded locally free resolution 
	by Lemma \ref{uniqueness}, 
	we find unique elements $P, P' \in \Matrix_{n+1}(k)$
	which make the following diagram commute:
	\begin{equation}\label{Diagram5}
		\xymatrix{
			0 \ar[r] &
			D \left( {\displaystyle \bigoplus_{i=0}^{n} \strshf_{\projsp^m}(-1)e_i} \right)[-1]
			 \ar[r]^(.6){\widetilde{M}} \ar[d]^{P'} &
			{\displaystyle \bigoplus_{i=0}^{n} \strshf_{\projsp^m}(-1)e_i} 
			\ar[r]^(.68){p} \ar[d]^{\transpose P}  &
			\mathcal{M} \ar[r] \ar[d]^{l} & 0 \\
			0 \ar[r] &
			D \left( {\displaystyle \bigoplus_{i=0}^{n} \strshf_{\projsp^m}(-1)e_i} \right)[-1]
			\ar[r]^(.6){\widetilde{M}} &
			{\displaystyle \bigoplus_{i=0}^{n} \strshf_{\projsp^m}(-1)e_i} 
			\ar[r]^(.68){p} &
			\mathcal{M} \ar[r] & 0.
		}
	\end{equation}
		
	This proves $L_0$ is embedded into $\Matrix_{n+1}(k) \times \Matrix_{n+1}(k)$
	as a $k$-subalgebra. The image satisfies $\transpose P M = M P'.$
	Conversely, if we take $(P, P')$ satisfying $\transpose P M = M P',$
	we obtain a morphism $l$ which makes the diagram \Equref{Diagram5} commute.
	Thus the image of $L_0$ consists of the set \Equref{endom1}.
	Transposing $\transpose P M = M P',$ we have
	$\transpose P' M = M P,$ and we see that $(P', P)$ is also an element of $L_0.$
	
	Let us take $l$ from $L.$
	Dualizing the diagram \Equref{Diagram5}, we obtain
	\[
		\xymatrix{
			0 \ar[r] &
			D \left( {\displaystyle \bigoplus_{i=0}^{n} \strshf_{\projsp^m}(-1)e_i} \right)[-1]
			 \ar[r]^(.6){\widetilde{M}}  &
			{\displaystyle \bigoplus_{i=0}^{n} \strshf_{\projsp^m}(-1)e_i} 
			\ar[r]^(.65){q}  &
			D\mathcal{M} \ar[r]  & 0 \\
			0 \ar[r] &
			D \left( {\displaystyle \bigoplus_{i=0}^{n} \strshf_{\projsp^m}(-1)e_i} \right)[-1]
			\ar[r]^(.6){\widetilde{M}}  \ar[u]^{P}&
			{\displaystyle \bigoplus_{i=0}^{n} \strshf_{\projsp^m}(-1)e_i} 
			\ar[r]^(.65){q} \ar[u]^{\transpose P'}&
			D\mathcal{M} \ar[r] \ar[u]^{\transpose l }& 0.
		}
	\]
	
	Now we would like to show
	\(
		P' = P.
	\)
	Since $q = \lambda \circ p,$ we have
	\begin{align*}
		\transpose l  \circ q &= 
		\transpose l  \circ \lambda \circ p 
		= \lambda \circ l \circ p = \lambda \circ p \circ \transpose P.
	\end{align*}
	Meanwhile, we have 
	\begin{align*}
		\transpose l  \circ q &=
		q \circ \transpose P'
		=  \lambda \circ p \circ \transpose P'.
	\end{align*}
	Since $\mathcal{M}(1)$ is generated by its global sections, we have 
	\(
		P' = P.
	\)
	Thus $L$ is contained in $L^\sigma_0.$
	Conversely, if we take $(P, P)$ corresponding to an element of $L_0^\sigma,$
	it is easy to see that the corresponding endomorphism $l$ is in $L.$ \qed
\end{prf}

Next, we study when $(\mathcal{M}, \lambda)$ is equivalent to $(\mathcal{M}, l\lambda)$ with respect to the equivalence relation $\sim_1$ (resp.\ $\sim_2$).

We treat the equivalence relation $\sim_1$ first.
Assume that $(\mathcal{M}, \lambda) \sim_1 (\mathcal{M}, l\lambda)$ 
for an element $l \in L.$
By definition, there exist an automorphism $\rho \colon
\mathcal{M} \longrightarrow \mathcal{M}$ and a constant $a \in k^\times$
such that 
\[
	l\lambda = a\transpose \rho  \circ \lambda \circ \rho.
\] 
Let $(P, P')$ be the element corresponding to 
$\rho \in L_0^\times,$
and $(R, R)$ the element corresponding to $l \in L.$  
Then $\sigma((P, P')) = (P', P)$ corresponds to 
$\sigma(\rho) := \lambda^{-1} \circ \transpose \rho  \circ \lambda.$ 
Hence we have
\begin{equation}\label{equality0}
	l = a \sigma(\rho) \rho \quad \Leftrightarrow \quad R = aPP'.
\end{equation}
Conversely, if we can write $l = a \sigma(\rho) \rho \in L$ 
for some $\rho \in L_0,$ we can go backward and obtain 
$(\mathcal{M}, \lambda) \sim_1 (\mathcal{M}, l\lambda).$

Next we treat the equivalence relation $\sim_2.$
Assume that $(\mathcal{M}, \lambda) \sim_2 (\mathcal{M}, l\lambda)$
for an element $l \in L.$
Take any ordered $k$-basis $s$ of $\Cohomology{0}{\projsp^m, \mathcal{M}(1)}.$
Put $M := \psi_c([(\mathcal{M}, \lambda, s)]).$
Then there exist an element $A \in \GenLin_{m+1}(k)$
and an isomorphism $\rho \colon \mathcal{M} \longrightarrow \nu_A^*\mathcal{M}$ 
such that
\begin{equation}\label{equality1}
	l\lambda = \transpose \rho  \circ \nu_A^*\lambda \circ \rho.
\end{equation}
There exists $(P, P') \in \Matrix_{n+1}(k) \times \Matrix_{n+1}(k)$ such that 
\[
	\xymatrix{
		0 \ar[r] &
		D\left( {\displaystyle \bigoplus_{i=0}^{n}	\strshf_{\projsp^m}(-1) e_i} \right)[-1] 
		 \ar[r]^(.6){M} \ar[d]_{P'} &
		{\displaystyle \bigoplus_{i=0}^{n} \strshf_{\projsp^m}(-1)e_i} 
		\ar[r]^(.7){p} \ar[d]_{\transpose P}  &
		\mathcal{M} \ar[r] \ar[d]^{\vertisom}_{\rho} & 0 \\
		0 \ar[r] &
		D \left( {\displaystyle \bigoplus_{i=0}^{n} \strshf_{\projsp^m}(-1) e_i} \right)[-1]
		\ar[r]^(.6){\widetilde{M} \cdot A} &
		{\displaystyle \bigoplus_{i=0}^{n} \strshf_{\projsp^m}(-1)e_i} 
		\ar[r]^(.65){\nu_A^* p} &
		\nu_A^* \mathcal{M} \ar[r] & 0
	}
\]
commutes.
Take $(R, R)$ corresponding to $l \in L.$ 
Then \Equref{equality1} means
\begin{equation}\label{equality2}
	M R = \transpose P M P \cdot A.
\end{equation}
We see that the converse holds.
For $(R, R)$ corresponding to an element $l \in L$ such that there exist
$(P, P') \in L_0, A \in \GenLin_{m+1}(k)$ satisfying \Equref{equality2},
then by going backward, we see that 
$(\mathcal{M}, \lambda) \sim_2 (\mathcal{M}, l\lambda).$
Thus we have proved the following proposition.

\begin{prop}\label{fiber}
	Let $M \in \NVLoci$ be an element, 
	and we write $\phi_c(M) = [(\mathcal{M}, \lambda, s)].$
	\begin{enumerate}[(i)]
		\item The $k$-algebra $L_0 = \Endomorphism{\projsp^m}{\mathcal{M}}$ acts on  
		$\Homomorphism{\projsp^m}{\mathcal{M}, 
		D\mathcal{M}}$ simply transitively in two ways.
		\item The $k$-subspace $L \subset L_0$ defined by the condition
		\(
			\transpose l \circ \lambda = \lambda \circ l
		\)
		is bijective to
		the $k$-subspace of symmetric homomorphisms of $\Homomorphism{\projsp^m}
		{\mathcal{M}, D{\mathcal{M}}}$ via
		\(
			l \mapsto \lambda \circ l.
		\)
		\item We define the subset $L_1$ of $L \cap L_0^\times$ consisting of $l$ such that 
		there exist $a \in k^\times, \rho \in L_0$ satisfying \Equref{equality0}.
		Then $L_1$ is bijective to the set of $\lambda'$ satisfying
		\(
			(\mathcal{M}, \lambda) \sim_1 (\mathcal{M}, \lambda').
		\)
		\item We define the subset $L_2$ of $L \cap L_0^\times$ consisting of $l$ such that
		there exist $A \in \GenLin_{m+1}(k)$ and 
		$\rho \colon \mathcal{M} \isomarrow \nu_A^*\mathcal{M}$ 
		satisfying \Equref{equality2}. Then $L_2$ is bijective to the set of $\lambda'$ satisfying
		\(
			(\mathcal{M}, \lambda) \sim_2 (\mathcal{M}, \lambda').
		\)
	\end{enumerate}
\end{prop}
\begin{rmk}
	The $k$-subspace $L \subset L_0$ depends on the choice of $\lambda.$ 
\end{rmk}

Let us assume two more conditions that $L_0$ is commutative and $\sigma = \identity.$
Under these assumptions, we have
\begin{align*}
	L_0 &= L,\\
	L_1 &= k^\times L^{\times 2} 
	= \left\{ ab^2 \relmid a \in k^\times, b \in L^\times \right\}.
\end{align*}
Hence we obtain the following corollary:

\begin{cor}\label{fiber2}
	Let $M \in \NVLoci$ be an element and  
	we write $\phi_c(M) = [(\mathcal{M}, \lambda, s)].$
	Moreover, we assume that $L_0$ is commutative and $\sigma = \identity.$
	Then the following statements hold.
	\begin{enumerate}[(i)]
		\item The set of symmetric quasi-isomorphisms from $\mathcal{M}$ 
		to $D\mathcal{M}$ has a simply transitive action of $L_0^\times.$ 
		\item The set of equivalence classes of symmetric quasi-isomorphisms from $\mathcal{M}$ 
		to $D\mathcal{M}$ with respect to the equivalence relation $\sim_1$
		has a simply transitive action of $L_0^\times /
		k^\times L_0^{\times 2}.$
	\end{enumerate}
\end{cor}

\section{Theta characteristics and the proof of Theorem \ref{main-theorem}}\label{Curves}
In this section, we discuss some properties of theta characteristics on 
geometrically reduced hypersurfaces and
prove a bijection which directly relates to the main theorem of this paper
(see Corollary \ref{bijection4}).
As before, 
we fix a quasi-isomorphism of complexes of coherent $\strshf_{\projsp^m}$-modules 
\[
	c \colon \omega_{\projsp^m} \isomarrow \strshf_{\projsp^m}(-m-1)[m].
\]

In Subsection \ref{bij2}, we state a bijection considering the multiplicity on
each irreducible component of the support of coherent $\strshf_{\projsp^m}$-modules. 
In Subsection \ref{TCs}, we introduce the notion of theta characteristics
on hypersurfaces, and give the desired bijection to show Theorem \ref{main-theorem}.
In Subsection \ref{main}, we prove Theorem \ref{main-theorem} and
Corollary \ref{Beauville-theorem} (see Corollary \ref{BeauTheo}).
In Subsection \ref{exaTC}, we give some examples of theta characteristics on plane curves.

\subsection{A bijection on $\strshf_{\projsp^m}$-modules 
with given supports and multiplicities}\label{bij2}
First, we recall some elementary facts from intersection theory
(\cite{Fulton-Intersection}).

Recall that, for a coherent $\strshf_{\projsp^m}$-module $\mathcal{F},$
there is a polynomial $P_{\mathcal{F}}(t)$ of degree less than or equal to 
$\dim \Supp(\mathcal{F})$ satisfying
\begin{align*}
	P_{\mathcal{F}}(t) &= \chi(\mathcal{F}(t)) \\
	&:= \sum_{i=0}^m (-1)^i \dim \Cohomology{i}{\projsp^m, \mathcal{F}(t)}.
\end{align*}
The polynomial $P_{\mathcal{F}}(t)$ is called the {\it Hilbert polynomial} of $\mathcal{F}.$

	\begin{lem}[{\cite[Example 2.5.2]{Fulton-Intersection}}]\label{H_0}
		Let $\mathcal{M}$ be a coherent $\strshf_{\projsp^m}$-module
		with $\dim \Supp(\mathcal{M}) = m-1.$
		The coefficient of $t^{m-1}$ in the Hilbert polynomial $P_\mathcal{M}(t)$
		is equal to
		\[
			\frac{1}{(m-1)!} \sum_{\eta \in \Generics(\Supp(\mathcal{M}))} \deg [\eta] \cdot
			\length_{\strshf_{\projsp^m,\eta}} (\mathcal{M}_\eta),
		\]
		where $[\eta]$ denotes the algebraic cycle on $\projsp^m$ corresponding to 
	the generic point $\eta$ of an irreducible component of $\Supp(\mathcal{M}).$
	\end{lem}
	\begin{prf}
		If we have a short exact sequence
		\[
			\xymatrix{
			0 \ar[r] & \mathcal{F}' \ar[r] & 
			\mathcal{F} \ar[r] & \mathcal{F}'' \ar[r] & 0 
			}
		\]
		and if the assertion of the lemma holds for $\mathcal{F}'$ and $\mathcal{F}'',$
		the assertion of the lemma for $\mathcal{F}$ 
		also holds by the additivity of the Hilbert polynomial and length of
		coherent $\strshf_{\projsp^m}$-modules.
		
		Let $\iota \colon Z \hookrightarrow \projsp^m$ be an irreducible component of
		$\Supp (\mathcal{M})$ of dimension $m-1,$
		and $\eta$ the generic point of $Z.$ 
		Then we have the following short exact sequence:
		\[
			\xymatrix{
				0 \ar[r] & \Kernel(\varphi) \ar[r] & 
				\mathcal{M} \ar[r]^(.45){\varphi} & \iota_*\iota^*\mathcal{M} \ar[r] & 0. 
			}
		\]
		The last term $\iota_*\iota^*\mathcal{M}$ is supported on $Z,$
		and we have
		\[
			\length_{\strshf_{\projsp^m}, \eta}(\Kernel(\varphi)) < 
			\length_{\strshf_{\projsp^m}, \eta}(\mathcal{M}).
		\]
		If $\dim \Supp(\mathcal{M}) < m-1,$ then the coefficient of $t^{m-1}$ is zero.
		Hence we may assume that $\mathcal{M}$ is supported on 
		an irreducible hypersurface $Z \subset \projsp^m.$
		
		The local ring $\strshf_{\projsp^m, \eta}$ is a discrete valuation ring. 
		Since $(\iota^*\mathcal{M})|_U$ is a free $\strshf_U$-module of finite rank
		for a non-empty Zariski open subscheme $U \subset Z.$ 
		we may assume that
		$\mathcal{M}$ is isomorphic to $\iota_*\strshf_Z.$
		In that case, the required equality is just the definition of the degree of $Z.$ \qed
	\end{prf}

\begin{lem}\label{intersections}
	Let $\mathcal{M}$ be a coherent $\strshf_{\projsp^m}$-module
	satisfying the following conditions:
	\begin{itemize}
		\item $\mathcal{M}$ is {\ArithCM}, and 
		\item $\mathcal{M}$ is pure of dimension $m-1,$ and
		\item $\Cohomology{0}{\projsp^m, \mathcal{M}} =0,$ and
		\item there exists a quasi-isomorphism of $\strshf_{\projsp^m}$-modules
		\[
			\lambda \colon \mathcal{M} \isomarrow
			\Right\Sheafhom{\projsp^m}{\mathcal{M}(2-m), \omega_{\projsp^m}[-m+1]}.
		\]
	\end{itemize} 
	Then we have
	\[
		\dim \Cohomology{0}{\projsp^m, \mathcal{M}(t)} = 
		\left( \sum_{\eta \in \Generics(\Supp(\mathcal{M}))} \deg [\eta] \cdot
		\length_{\strshf_{\projsp^m,\eta}} (\mathcal{M}_\eta) \right)
		\binom{t+m-2}{m-1}
	\]
	for $t \ge 0.$
\end{lem}

\begin{prf}
	By Proposition \ref{ACM-and-resolution} and Lemma \ref{Cohomology},
	the coherent $\strshf_{\projsp^m}$-module $\mathcal{M}$ admits 
	a short exact sequence of the following form:
	\[\xymatrix{
		0 \ar[r] & {\displaystyle \bigoplus_{i=0}^{r} \strshf_{\projsp^m}(-2)} \ar[r]
		& {\displaystyle \bigoplus_{i=0}^{r} \strshf_{\projsp^m}(-1)} \ar[r]
		& \mathcal{M} \ar[r] & 0
	}\]
	for some $r \ge 0.$
	From the long exact sequence of cohomology, we have
	$\Cohomology{i}{\projsp^m, \mathcal{M}(t)} = 0$ for $i \ge 1$ and $t \ge 0.$
	By the additivity of the Hilbert polynomials, we have
	\[
		P_{\mathcal{M}}(t) =
		(r+1) \left( P_{\strshf_{\projsp^m}}(t-1) - P_{\strshf_{\projsp^m}}(t-2) \right).
	\]
	
	The Hilbert polynomial $P_{\strshf_{\projsp^m}}(t)$ is known as
	\[
		P_{\strshf_{\projsp^m}}(t) = \binom{t + m}{m}.
	\]
	Hence we have
	\begin{align*}
		P_{\mathcal{M}}(t) &=
		(r+1) \left( \binom{t+m-1}{m} - \binom{t+m-2}{m} \right) \\
		&= (r+1) \binom{t+m-2}{m-1}.
	\end{align*}
	Comparing the coefficient by $t^{m-1}$ by Lemma \ref{H_0}, 
	we obtain the desired equality. \qed
\end{prf}

By Theorem \ref{bijection1} and  Lemma \ref{intersections},
we obtain the following proposition.

\begin{prop}\label{bijection2}
	Let $S_1, S_2, \dots, S_r$ be a collection of distinct 
	irreducible hypersurfaces in $\projsp^m$ over $k,$
	$\eta_1, \eta_2, \dots, \eta_r$ their generic points, 
	$F_1, F_2, \dots, F_r$ their defining equations and
	$n_1, n_2, \dots, n_r$ non-negative integers.
	Put $\mathbf{n}= (n_1,n_2, \dots, n_r)$ and 
	write $|\mathbf{n}| = \sum_i n_i \deg [\eta_i ] - 1.$ 
	Let $S := \bigcup_{i} S_i.$ 
	Then there exists a natural bijection between the following two sets.
	\begin{itemize}
		\item\label{matrices2}
		The set $U_{S, \mathbf{n}}$ of $(m+1)$-tuples of symmetric matrices 
		$M= (M_0, M_1, \dots, M_m)$ of size $(|\mathbf{n}|+1)$
		with entries in $k$ satisfying 
		\[
			\disc (M) = u F_1^{n_1}\cdot F_2^{n_2} \cdot \dots \cdot F_r^{n_r} 
		\] 
		for some $u \in k^\times.$
		\item\label{sheaves2}
		The set $V_{S,\mathbf{n}}$ of equivalence classes of triples 
		$(\mathcal{M}, \lambda, s),$
		where
		\begin{itemize}
			\item $\mathcal{M}$ is a coherent $\strshf_{\projsp^m}$-module
				satisfying the following conditions:
				\begin{itemize}
					\item $\mathcal{M}$ is {\ArithCM,} and
					\item $\mathcal{M}$ is pure of dimension $m-1,$ and
					\item $\Supp(\mathcal{M}) \subset S,$ and
					\item $\length_{\strshf_{\projsp^m, \eta_i}}( \mathcal{M}_{\eta_i}) = n_i$ 
					for any $1 \le i \le r,$ and
					\item $\Cohomology{0}{\projsp^m, \mathcal{M}}=0$ and
					$\dim \Cohomology{0}{\projsp^m, \mathcal{M}(1)} = |\mathbf{n}|+1.$
				\end{itemize}
			\item $\lambda$ is a symmetric quasi-isomorphism
				\[
					\lambda \colon \mathcal{M} \isomarrow
					 \Right\Sheafhom{\projsp^m}{\mathcal{M}(2-m), \omega_{\projsp^m}
					 [-m+1]}.
				\]
			\item $s= \{ s_0, s_1, \dots, s_{|\mathbf{n}|} \}$ is an ordered $k$-basis of
				$\Cohomology{0}{\projsp^m, \mathcal{M}(1)}.$
			\end{itemize}
		Here, two triples $(\mathcal{M}, \lambda, s), (\mathcal{M}', \lambda', s')$ are 
		said to be {\it equivalent} if
		there exists an isomorphism $\rho \colon \mathcal{M} \isomarrow \mathcal{M}'$ 
		of $\strshf_{\projsp^m}$-modules satisfying 
		\begin{itemize}
			\item 
			\(
				\transpose \rho  \circ \lambda' \circ \rho
			= \lambda, 
			\) and
			\item 
			\( 
				\rho(s_i) = s'_i
			\)
			for any $0 \le i \le |\mathbf{n}|.$
		\end{itemize}
	\end{itemize}
\end{prop}

\begin{prf}
	Let us take an element $M \in U_{S, \mathbf{n}}.$
	Then, as before, we have the following exact sequence
	\[
		\xymatrix{
			0 \ar[r] & {\displaystyle \bigoplus_{i=0}^{|\mathbf{n}|} 
			\strshf_{\projsp^m}(-2) }
			\ar[r]^(.5){M} &
			{\displaystyle \bigoplus_{i=0}^{|\mathbf{n}|} \strshf_{\projsp^m}(-1)} 
			\ar[r]^(.7){p} &
			\mathcal{M} \ar[r] & 0.
		}
	\]
	Localizing this sequence at the generic point $\eta_i$ of $S_i,$
	we have the short exact sequence of $\strshf_{\projsp^m, \eta_i}$-modules:
	\[
		\xymatrix{
			0 \ar[r] & {\displaystyle \bigoplus_{i=0}^{|\mathbf{n}|} 
			\strshf_{\projsp^m, \eta_i}(-2) }
			\ar[r]^(.5){M} &
			{\displaystyle \bigoplus_{i=0}^{|\mathbf{n}|} \strshf_{\projsp^m, \eta_i}(-1)} 
			\ar[r]^(.68){p} &
			\mathcal{M}_{\eta_i} \ar[r] & 0.
		}
	\]
	By \cite[Lemma A.2.6]{Fulton-Intersection},
	we have
	\begin{align*}
		\length_{\strshf_{\projsp^m, \eta_i}}(\mathcal{M}_{\eta_i})
		&= \length_{\strshf_{\projsp^m, \eta_i}}(\strshf_{\projsp^m, \eta_i}/\disc(M)) \\
		&= \mathrm{ord}_{F_i}(\disc(M)) \\
		&= n_i.
	\end{align*}
	Hence $\phi_c(M) \in V_{S, \mathbf{n}}.$ 
	
	Conversely, take a representative $(\mathcal{M}, \lambda, s)$ of an element of 
	$V_{S, \mathbf{n}},$ and we put $M=\psi_c([(\mathcal{M}, \lambda, s)]).$
	Again by \cite[Lemma A.2.6]{Fulton-Intersection}, 
	we have 
	\[
		\mathrm{ord}_{F_i}(\disc(M))=n_i.
	\]
	Since we have assumed that
	\[
		\dim \Cohomology{0}{\projsp^m, \mathcal{M}(1)} = |\mathbf{n}|+1,
	\]
	we have 
	\begin{align*}
		\deg \disc(M) &= |\mathbf{n}|+1 \\
		&=\sum_i n_i \deg[\eta_i] \\
		&= \sum_i n_i \deg F_i
	\end{align*}
	by Lemma \ref{intersections}. Hence we have
	\[
		\disc(M) = u F_1^{n_1} \cdot F_2^{n_2} \cdot \dots \cdot
		F_r^{n_r}
	\]
	for some $u \in k^\times.$ Thus $\psi_c([(\mathcal{M}, \lambda, s)]) \in 
	U_{S, \mathbf{n}}.$
\qed\end{prf}

\subsection{Theta characteristics on hypersurfaces}\label{TCs}
In this subsection, for a geometrically reduced hypersurface $S \subset \projsp^m,$
we introduce the notion of theta characteristics on $S$ following Mumford, Harris, Piontkowski,
Dolgachev.

\begin{dfn}[{\cite{Mumford-Theta}, \cite{Harris-Theta}, \cite{Piontkowski-Theta}, 
\cite[Definition 4.2.9]{Dolgachev-CAG}}]\label{TC}
	Let $\iota \colon S \hookrightarrow \projsp^m$ be 
	a geometrically reduced hypersurface over $k.$ 
	A {\it theta characteristic} on $S$ is a coherent $\strshf_S$-module $\mathcal{M}$ 
	such that $\iota_*\mathcal{M}$ is {\ArithCM}, 
	pure of dimension $m-1,$
	\[
		\length_{\strshf_{S, \eta}}(\mathcal{M}_\eta) =1
	\]
	for each generic point $\eta \in \Generics(S),$
	and there is a quasi-isomorphism of complexes of coherent $\strshf_S$-modules
	\[
		\lambda \colon 
		\mathcal{M} \isomarrow \Right\Sheafhom{S}{\mathcal{M}(2-m), \omega_S[-m+1]}.
	\]
	A theta characteristic $\mathcal{M}$ on $S$ is said to be {\it effective} 
	(resp.\ {\it non-effective})
	if $\Cohomology{0}{S, \mathcal{M}} \neq 0$ (resp.\ 
	$\Cohomology{0}{S, \mathcal{M}} =0$).
\end{dfn}

\begin{rmk}
	Let $k$ be an algebraically closed field.
	Note that, though smooth plane curves over $k$ have theta characteristics
	(see Example \ref{curvesTC}),
	general hypersurfaces of dimension $\ge 2$ over $k$ do not have theta characteristics
	(cf.\ \cite[Corollary 6.6]{Beauville-Det}, \cite[Example 4.2.23]{Dolgachev-CAG}). 
\end{rmk}

\begin{rmk}
	We can modify the definition of theta characteristics
	by the existence of an isomorphism 
	\[
		\lambda^{\mathrm{shf}} \colon \mathcal{M} 
		\isomarrow \Sheafhom{S}{\mathcal{M}(2-m), \omega_S^{\mathrm{shf}}}
	\]
	of $\strshf_S$-modules	instead of the quasi-isomorphism $\lambda,$
	where $\omega_S^{\mathrm{shf}}$ denotes the {\it dualizing sheaf} on $S.$
	Since $S$ is a hypersurface, its dualizing sheaf is a line bundle.
	Since $\iota_*\mathcal{M}$ is {\ArithCM} and pure of dimension $m-1,$ 
	$\iota_*\mathcal{M}$ has a minimal graded locally free resolution 
	of length one; namely, we have
	\[
		\xymatrix{
			0 \ar[r] & {\displaystyle \bigoplus_{i=0}^r \strshf_{\projsp^m}(e_i)} \ar[r] &
			{\displaystyle \bigoplus_{i=0}^r \strshf_{\projsp^m}(d_i)}
			\ar[r] & \iota_* \mathcal{M} \ar[r] & 0
		}
	\]
	for some integers $r \ge 0$ and $d_i, e_i \in \Integer$ 
	(\cite[Proposition 1.2]{Beauville-Det}). 
	Since $\iota_*\mathcal{M}$ is a torsion sheaf on $\projsp^m,$
	we have
	\[
		\Sheafhom{\projsp^m}{\iota_*\mathcal{M}(2-m),
		\omega_{\projsp^m}^{\mathrm{shf}}} = 0.
	\]
	On the other hand, we have
	\[
		\Sheafhom{\projsp^m}{\bigoplus_{i=0}^r \strshf_{\projsp^m}(2-m+j),
		\omega_{\projsp^m}^{\mathrm{shf}}} \cong 
		\bigoplus_{i=0}^r \strshf_{\projsp^m}(-3-j).
	\]
	Hence applying the functor $\Sheafhom{\projsp^m}{*(2-m),
	 \omega_{\projsp^m}^{\mathrm{shf}}}$ to the above short exact sequence,
	we obtain 
	\[
		\xymatrix{
			0 \ar[r] & {\displaystyle \bigoplus_{i=0}^r \strshf_{\projsp^m}(e_i)} \ar[r] &
			{\displaystyle \bigoplus_{i=0}^r \strshf_{\projsp^m}(d_i)}
			\ar[r] & \Sheafext{1}{\projsp^m}{\iota_* \mathcal{M}(2-m),
			\omega_{\projsp^m}^{\mathrm{shf}}} \ar[r] & 0
		}
	\]
	and
	\[
		\Sheafext{i}{\projsp^m}
		{\iota_*\mathcal{M}(2-m), \omega_{\projsp^m}^{\mathrm{shf}}} = 0 
		\quad (i \neq 1, j \in \Integer).
	\]
	These show 
	\[
		\Right\Sheafhom{\projsp^m}{\iota_*\mathcal{M}(2-m), 
		\omega_{\projsp^m}[-m+1]} \cong 
		\Sheafext{1}{\projsp^m}{\iota_* \mathcal{M}(2-m),
		\omega_{\projsp^m}^{\mathrm{shf}}}
	\]
	and by Grothendieck duality
	\[
		\Right\Sheafhom{S}{\mathcal{M}(2-m), 
		\omega_{S}[-m+1]} \cong 
		\Sheafhom{S}{\mathcal{M}(2-m),
		\omega_{S}^{\mathrm{shf}}}.
	\]
	Hence $\Right\Sheafhom{S}{\mathcal{M}(2-m), \omega_S[-m+1]}$ 
	is quasi-isomorphic to $\Sheafhom{S}{\mathcal{M}(2-m), \omega_S^{\mathrm{shf}}}.$
	This shows that, if an isomorphism $\lambda^{\mathrm{shf}}$ exists,
	it can be interpreted as the quasi-isomorphism $\lambda.$
	In particular, when $m=2$ and $S \subset \projsp^2$ is a geometrically reduced plane curve,
	the definition of theta characteristics in Definition \ref{TC},
	coincides with the Piontkowski's definition of 
	theta characteristics on reduced singular curves (\cite{Piontkowski-Theta}). 
\end{rmk}

\begin{rmk}\label{symmetricity2}
	If a quasi-isomorphism $\lambda$ exists, it is automatically symmetric.
	In fact, it is enough to show that any morphism 
	\[
		v \colon \mathcal{M} \to
		\Right\Sheafhom{S}{\mathcal{M}(2-m), \omega_{S}[-m+1]}
	\]
	is symmetric. To prove it, consider the morphism  
	\[
		h_v = 	\transpose v  \circ \canonical - v.
	\]
	Since $\mathcal{M}$ is locally free of rank one on 
	a Zariski open dense subscheme $U \subset S,$ 
	the morphism $h_v$ is zero on $U.$
	So the $\strshf_S$-submodule $\Image (h_v)$ of $\Right\Sheafhom{S}
	{\mathcal{M}(2-m), \omega_S[-m+1]} \cong \mathcal{M}$ is 
	supported on a subscheme of $S$ whose dimension is strictly less than
	$\dim S = m-1.$
	Since $\mathcal{M}$ is pure of dimension $m-1,$ 
	we conclude $\Image(h_v) = 0$
	and $v$ is symmetric. (See also Lemma \ref{End2}.)
\end{rmk}

Let $\mathscr{I}_S$ be the ideal sheaf of $\strshf_{\projsp^m}$ defining 
$\iota \colon S \hookrightarrow \projsp^m.$
The push-forward functor $\iota_*$ gives a categorical equivalence
between the category of coherent $\strshf_{\projsp^m}$-modules 
which are annihilated by $\mathscr{I}_S$ 
and the category of coherent $\strshf_{S}$-modules.
Hence, combining this with Corollary \ref{Bij-cor1} and Proposition \ref{bijection2}
for $n_1 = n_2 = \dots = n_r =1,$ we obtain the following corollary:

\begin{cor}\label{bijection4}
	Let $S \subset \projsp^m$ be a geometrically reduced hypersurface over $k.$ 
	Then there exists a natural bijection between the following two sets.
	\begin{itemize}
		\item\label{matrices4}
		The set of $(k^\times I_{m+1}) \times \GenLin_{n+1}(k)$-orbits of
		$(m+1)$-tuples of symmetric matrices 
		$M= (M_0, M_1, \dots, M_m)$ of size $n+1$ with entries in $k$ such that
		the equation $(\disc(M)=0)$ defines $S.$
		\item\label{sheaves}
		The set $\ThetaChar_{m+1, n+1}(k)_S^\sim$ of equivalence classes of pairs 
		$(\mathcal{M}, \lambda)$ with respect to the equivalence relation $\sim_1,$ where
		\begin{itemize}
			\item $\mathcal{M}$ is a non-effective theta characteristic on $S,$ and
			\item  $\lambda$ is a quasi-isomorphism
				\[
					\lambda \colon \mathcal{M} \isomarrow
					\Right\Sheafhom{S}{\mathcal{M}(2-m), \omega_{S}[-m+1]}.
				\]
			\end{itemize}
	\end{itemize}
\end{cor}
\begin{prf}
	Recall that any quasi-isomorphism $\lambda$
	is symmetric by Remark \ref{symmetricity2}.
	Hence we only have to show that the cokernel of $M \in \RedLoci$ can be lifted to
	a coherent $\strshf_S$-module.

	Take $M \in \RedLoci$ and put $\phi_c(M) = [(\mathcal{M}, \lambda, s)].$
	It is enough to show that 
	\[
		\disc(M)\mathcal{M} = 0.
	\]
	
	Take a minimal graded locally free resolution of $\mathcal{M},$ 
	and consider the following commutative diagram:
	\[
		\xymatrix{
			0 \ar[r] &
			{\displaystyle \bigoplus_{i=0}^{n} \strshf_{\projsp^m}(-n-3)}
			 \ar[r]^(.5){M} \ar[d]^{\times \disc(M)} &
			{\displaystyle \bigoplus_{i=0}^{n} \strshf_{\projsp^m}(-n-2)} 
			\ar[r]^(.58){q} \ar[d]^{\times \disc(M)} &
			\mathcal{M}(-n-1) \ar[r] \ar[d]^{\times \disc(M)} & 0 \\
			0 \ar[r] &
			{\displaystyle \bigoplus_{i=0}^{n} \strshf_{\projsp^m}(-2)}
			\ar[r]^(.5){M}  &
			{\displaystyle \bigoplus_{i=0}^{n} \strshf_{\projsp^m}(-1)} 
			\ar[r]^(.63){q}&
			\mathcal{M} \ar[r] & 0.
		}
	\]
	To prove the condition 
	\(
		\disc(M)\mathcal{M} = 0,
	\)
	it suffices to show
	\[
		\disc(M)\bigoplus_{i=0}^{n} \strshf_{\projsp^m}(-n-2) \subset
		M\bigoplus_{i=0}^{n} \strshf_{\projsp^m}(-2).
	\]
	Take the adjugate matrix $\mathrm{adj}({M})$ of $M,$
	then we have
	\[
		M\bigoplus_{i=0}^{n} \strshf_{\projsp^m}(-2) \supset
		M\mathrm{adj}({M})\bigoplus_{i=0}^{n} \strshf_{\projsp^m}(-n-2) 
		= \disc(M)\bigoplus_{i=0}^{n} \strshf_{\projsp^m}(-n-2).
	\]
	This concludes the proof. \qed
\end{prf}

\begin{rmk}
As a consequence, we obtain a bijection concerning symmetric determinantal representations
of a geometrically reduced hypersurface $S$ of degree $n+1$ over $k.$ 
Here, if $F$ is a defining equation of $S,$ 
a {\it symmetric determinantal representation} of $S$ is 
an $(m+1)$-tuple of symmetric matrices $M$ of size $n+1$ with
\[
	\disc(M) = uF
\]
for some $u \in k^\times.$ Two symmetric determinantal representations $M, M'$ are said to be
{\it equivalent} if there exist $P \in \GenLin_{n+1}(k)$ and $a \in k^\times$ such that
\[
	M' = a \transpose P M P.
\]
By Corollary \ref{bijection4}, the equivalence classes of 
symmetric determinantal representations of $S$ are in bijection with
the elements in $\ThetaChar_{m+1, n+1}(k)_S^\sim.$
We investigate symmetric determinantal representations of
smooth plane curves over global fields in \cite{Ishitsuka-Ito-Det1},\ 
\cite{Ishitsuka-Ito-Det2},\ \cite{Ishitsuka-Ito-Det3}.

\end{rmk}

\subsection{Proof of Theorem \ref{main-theorem} and 
its corollary}\label{main}
We define the set $\ThetaChar_{m+1, n+1}(k)^\sim$ as
\[
	\ThetaChar_{m+1, n+1}(k)^\sim := \bigcup_{S} \ThetaChar_{m+1, n+1}(k)_S^\sim,
\]
where $S$ runs over all geometrically reduced hypersurfaces in $\projsp^m$ defined over $k.$

We also define the set $\ThetaChar_{m+1, n+1}(k)$ as
the set of the equivalence classes of pairs $(S, \mathcal{M}),$
where two pairs $(S, \mathcal{M}), (S', \mathcal{M}')$ are said to be {\it equivalent} if
$S=S'$ and $\mathcal{M}, \mathcal{M}'$ are isomorphic as $\strshf_S$-modules. 
Since any quasi-isomorphism
\[
	\lambda \colon \mathcal{M} \isomarrow 
	\Right\Sheafhom{S}{\mathcal{M}(2-m), \omega_S[-m+1]}
\]
is symmetric (see Remark \ref{symmetricity2}), 
we can consider the set $\ThetaChar_{m+1, n+1}(k)$ as
the set obtained by forgetting the data of quasi-isomorphisms $\lambda.$
Hence we have a natural surjection
\[
	\ThetaChar_{m+1,n+1}(k)^\sim \longrightarrow \ThetaChar_{m+1, n+1}(k).
\]

We also define the map
\[
	\Phi_{m+1, n+1} \colon \RedLoci \longrightarrow \ThetaChar_{m+1, n+1}(k)
\] 
so that $\Phi_{m+1, n+1}(M)$
is the equivalence class $[(S, \mathcal{M})] \in \ThetaChar_{m+1, n+1}(k)$ of 
a pair $(S, \mathcal{M})$ satisfying
$\phi_c(M) = [(\mathcal{M}, \lambda,s)]$
and $\Supp (\mathcal{M}) =S.$ 
The surjectivity of $\Phi_{m+1, n+1}$ follows from Corollary \ref{bijection4}, 
so we only have to study the fibers of $\Phi_{m+1, n+1}.$ In other words,
we shall investigate the equivalence classes of symmetric quasi-isomorphisms 
$\lambda$ with respect to $\sim_1.$ 

%

\begin{lem}\label{End2}
	Take any geometrically reduced hypersurface $S \subset \projsp^m$ over $k,$ 
	and let $\mathcal{M}$ be a coherent
	$\strshf_S$-module satisfying the conditions in Proposition \ref{bijection2}.
	The endomorphism sheaf $\Sheafend{S}{\mathcal{M}}$ is
	embedded into $\prod_{\eta \in \Generics(S)} 
	\Sheafend{S}{i_{\eta, *}\mathcal{M}_\eta}.$
	In particular, if $\mathcal{M}$ is a theta characteristic on $S,$ 
	the $\strshf_S$-algebra $\Sheafend{S}{\mathcal{M}}$ is commutative
	and the $k$-algebra of its global sections 
	$L_0 = \Endomorphism{S}{\mathcal{M}} = L$ is an {\'e}tale $k$-algebra
	of finite degree.
\end{lem}
\begin{prf}
	Take any non-zero element $f$ in $\Sheafend{S}{\mathcal{M}}(U) = 
	\Endomorphism{U}{i_U^*\mathcal{M}}$ 
	for an open subscheme $i_U \colon U \hookrightarrow S.$
	Then $f$ also gives an endomorphism 
	$f_\eta$ of $i_{\eta, *}\mathcal{M}_\eta$ for each $\eta \in \Generics(U).$
	
	Since $S$ is reduced, if $f_\eta$ is the zero endomorphism for each $\eta \in \Generics(U),$
	the image $\Image(f) \subset i_{U}^{*}\mathcal{M}$ is supported on 
	a subscheme of $U$ of dimension less than $m-1.$
	Since $i_U^*\mathcal{M}$ is pure of dimension $m-1,$ 
	we have $f = 0.$
	
	If $\mathcal{M}$ is a theta characteristic on $S,$
	$\Sheafend{S}{i_{\eta, *}\mathcal{M}_\eta}$ is
	isomorphic to $i_{\eta, *}\strshf_{S, \eta}$ since $\mathcal{M}_\eta$ is a free 
	$\strshf_{S, \eta}$-module of rank one. 
	In particular, $\Sheafend{S}{i_{\eta, *}\mathcal{M}_\eta}$ is commutative
	and moreover isomorphic to the coordinate ring of $\eta,$ which is a field.
	Hence the product of the rings of global endomorphisms $\prod_{\eta \in \Generics(S)} 
	\Endomorphism{S}{i_{\eta, *}\mathcal{M}_\eta}$ is commutative,
	and the subring $L_0$ is also commutative.
	
	On the other hand, for any element $l \in L_0,$ the difference $\sigma(l) - l$
	is zero at each generic point $\eta \in \Generics(S).$
	Hence the dimension of the support of $\Image(\sigma(l)-l)$ is less than $m-1.$ 
	By the purity of $\mathcal{M},$ we have $\sigma(l) = l.$ 
	This shows $L_0 = L.$
	
	Finally we shall show that $L_0$ is an {\'e}tale $k$-algebra.
	By descent theory, we may assume that $k$ is algebraically closed.
	The algebra $L_0$ is a finite dimensional $k$-subalgebra of
	the product of fields 
	$\prod_{\eta \in \Generics(S)} \Endomorphism{S}{i_{\eta, *}\mathcal{M}_\eta}.$
	Hence $L_0$ has no nonzero nilpotent elements, and $L_0$ is an {\'e}tale $k$-algebra.
	\qed
\end{prf}
\begin{prf}[Proof of Theorem \ref{main-theorem}]
We fix a geometrically reduced hypersurface $S$ over $k,$ and
a theta characteristic $\mathcal{M}$ on $S.$
By Corollary \ref{fiber2} and Lemma \ref{End2},
the set of symmetric quasi-isomorphisms 
\[
	\mathcal{M} \isomarrow 
	\Right\Sheafhom{S}{\mathcal{M}(2-m), \omega_S[-m+1]}
\]
admits a transitive action of $L^\times,$
and we have $L_1 = k^\times L^{\times 2}.$
Since the fiber $\Phi_{m+1, n+1}^{-1}([(S, \mathcal{M})])$ is 
the set of equivalence classes with respect to $\sim_1,$
it admits a simply transitive action of $L^\times / k^\times L^{\times 2}.$
This completes the proof of Theorem \ref{main-theorem}.
\qed\end{prf}

\begin{cor}\label{BeauTheo}
	The fiber $\Phi_{m+1, n+1}^{-1}([(S, \mathcal{M})])$ is a singleton 
	if {\it at least one} of the following conditions is satisfied:
   \begin{itemize}
  		\item the base field $k$ is separably closed of characteristic different from two, or
		\item the base field $k$ is perfect of characteristic two, or
		\item the hypersurface $S \subset \projsp^m$ is geometrically integral.
   \end{itemize}
\end{cor}

\begin{prf}
	Recall that $L = L_0$ is an {\'e}tale $k$-algebra. 
	It is a product of finite separable extensions of $k.$
	If the first or second condition is satisfied, 
	we have $L^\times = L^{\times 2}$ and 
	the group $L^\times / k^\times L^{\times 2}$ is trivial.
	Since the fiber $\Phi_{m+1, n+1}^{-1}([(S, \mathcal{M})])$ has 
	a simply transitive action of the trivial group, it is a singleton.
	We assume that the third condition is satisfied.
	Then $L \otimes_k \overline{k}$ can be embedded into a field, 
	where $\overline{k}$ is an algebraic closure of $k.$
	Hence we conclude that $L \otimes_k \overline{k}$ is an integral domain.
	Hence we have $L=k.$
	The group $L^\times / k^\times L^{\times 2}$ is trivial, and
	the fiber $\Phi_{m+1, n+1}^{-1}([(S, \mathcal{M})])$ is a singleton.\qed
\end{prf}

\begin{rmk}
The endomorphism sheaf $\Sheafend{S}{\mathcal{M}}$ is called
the {\it global invariant} of a theta characteristic $\mathcal{M}$ 
(\cite[Definition 4.2.9]{Dolgachev-CAG}).
\end{rmk}

\subsection{Theta characteristics on plane curves}\label{exaTC}
Classical examples of theta characteristics are those on plane curves.
We give some examples. (See \cite{Mumford-Theta}, \cite{Harris-Theta}, 
\cite{Piontkowski-Theta} for details.)

\begin{exa}\label{curvesTC}
	If $C$ is a smooth plane curve of genus $g,$ 
	the definition of theta characteristics in Definition \ref{TC}
	coincides with the usual definition of theta characteristics on projective smooth curves
	due to Mumford (\cite{Mumford-Theta}).
	If $k$ is an algebraically closed field of characteristic different from two,
	there exist $2^{2g}$ theta characteristics on $C,$ 
	and the number of non-effective theta characteristics is less than or equal to
	$2^{g-1}(2^g+1).$
	If $k$ is algebraically closed of characteristic zero,
	there exists at least one non-effective theta characteristic on $C$ 
	(\cite[Remark 4.4]{Beauville-Det}).
\end{exa}

\begin{exa}[{\cite{Harris-Theta}, \cite{Piontkowski-Theta}}]\label{push-theta}
	We put $m=2.$ Take a geometrically reduced plane curve $C \subset \projsp^2.$
	For a partial normalization $\pi \colon N \to C,$
	the push-forward of a theta characteristic $\mathcal{L}$ on $N$ is
	a theta characteristic on $C.$ 
	Let $\mathcal{L}$ be a theta characteristic on $N,$ and
	\[
		\lambda \colon \mathcal{L} \isomarrow 
		\Right\Sheafhom{N}{\mathcal{L}, \omega_N[-1]}
	\]
	a quasi-isomorphism. Then, 
	by Grothendieck duality, we have the following quasi-isomorphisms on $C$
	\[
		\lambda_{\pi_*\mathcal{L}} \colon 
		\pi_*\mathcal{L} \underset{\pi_*\lambda}{\isomarrow}
		\pi_* \Right\Sheafhom{N}{\mathcal{L}, \omega_N[-1]} 
		\underset{\GDual_{\pi}}{\isomarrow}
		\Right\Sheafhom{C}{\pi_* \mathcal{L}, \omega_C[-1]}.
	\]
%
\end{exa}

\begin{exa}\label{infinite}
	The following example shows that the fiber of $\Phi_{m+1, n+1}$
	can have infinitely many elements. We put $m=2$ and $n=1.$
	Consider the union of two lines $S = (xy=0) \subset \projsp^2$ defined over $k=\Rational.$
	We can find a non-effective theta characteristic $\mathcal{M}$ on $S$
	which is the push-forward of a line bundle on the normalization of $S.$
	Then we have $L \cong \Rational \times \Rational.$ 
	In this case, 
	$L^\times / k^\times L^{\times 2} \cong \Rational^\times /\Rational^{\times 2}$ 
	is an infinite group, and the fiber $\Phi_{3, 2}^{-1}([(S, \mathcal{M})])$
	has infinitely many elements.
\end{exa}

\section{Projective automorphism groups of complete intersections of quadrics}\label{ProjAut}
In this section, we give an application of our methods
to the projective automorphism groups of complete intersections of quadrics
in the projective space. 

We fix a field $k$ of {\it characteristic different from two} in this section.
We also fix integers $n > m \ge 2.$
As before, we fix a quasi-isomorphism of complexes of 
coherent $\strshf_{\projsp^m}$-modules
\[
	c \colon \omega_{\projsp^m} \isomarrow \strshf_{\projsp^m}(-m-1)[m].
\]

\subsection{Preliminaries on quadrics}\label{Prel-of-qdrcs}
First, let us recall some notations and basic facts on quadrics 
(cf.\ \cite[Chapter 6]{Griffiths-Harris-AG}).

A {\it quadric} $Q_0$ in $\projsp^n$ means a hypersurface of degree two in $\projsp^n.$
We use the lower-case letter $q_0$ to denote one of the quadratic forms in $(n+1)$-variables
defining $Q_0.$
Let $M_{Q_0}$ be the Gram matrix of the quadratic form $q_0.$ 
Explicitly, we have
\[
	\transpose{x}M_{Q_0}x = q_0(x)
\] 
and 
\[
	Q_0 = (q_0(x)=0) \subset \projsp^n.
\]
Here, we identify points on $\projsp^n$
with $(n+1)$-dimensional column vectors.
Since the characteristic of $k$ is different from two,
quadrics in $\projsp^n$ are in one-to-one correspondence with points in the projective space 
$\projsp(\Sym_2 k^{n+1})$ of one-dimensional subspaces in 
the $k$-vector space of symmetric matrices of size $n+1.$

The {\it corank} of a quadric $Q_0$, denoted as $\corank(Q_0)$, is the corank of $M_{Q_0}$
 defined by 
 \[
 	\corank(Q_0) := \corank(M_{Q_0}) := n+1-\rank(M_{Q_0}).
 \]
It is classically known that $Q_0$ is a non-singular variety if and only if 
$\corank (Q_0) =0. $ If $Q_0$ is singular, the singular locus $\Singular (Q_0)$ is 
a $(\corank (Q_0) -1)$-dimensional linear subvariety of $\projsp^n.$
Actually, $\Singular(Q_0)$ is defined by the kernel of $M_{Q_0}$ as a $k$-linear map.
The locus 
\[
	\Delta \subset \projsp(\Sym_2 k^{n+1})
\]
of singular quadrics is defined by $(\det (M_{Q_0}) = 0),$ so $\Delta$ is a hypersurface
in $\projsp(\Sym_2 k^{n+1}).$

Let $\Pi_Q$ be an $m$-dimensional linear subvariety of $\projsp(\Sym_2k^{n+1}).$
Let $Q = (Q_0, Q_1, \dots, Q_m)$ be an $(m+1)$-tuple of quadrics in $\projsp^n$ 
spanning the linear subvariety $\Pi_Q.$ We write 
\[
	\Pi_Q = \langle Q_0, Q_1, \dots, Q_m \rangle \cong \projsp^m.
\]
The linear subvariety $\Pi_Q$ is called the {\it linear system of quadrics} 
generated by $Q.$ We denote the locus of singular quadrics in $\Pi_Q$ by 
\[
	\Delta_Q := \Delta \cap \Pi_Q.
\]
It is a hypersurface in $\Pi_Q \cong \projsp^m$ defined by the vanishing of the determinant
of Gram matrices.
We write the base locus of $\Pi_Q$ as 
\[
	X_Q := Q_0 \cap Q_1 \cap \dots \cap Q_m \subset \projsp^n.
\]

Let $M_Q \in W = k^{m+1} \otimes \Sym_2 k^{n+1}$ be
the symmetric matrix corresponding to an $(m+1)$-tuple of Gram matrices
$ (M_{Q_0}, M_{Q_1}, \dots, M_{Q_m}) .$
If $X_Q \subset \projsp^n$ is a complete intersection of 
$m+1$ quadrics $Q_0, Q_1, \dots,  Q_m,$
we see that $\codimension_{\projsp^n} X_Q = m+1$ and
the set of quadrics containing $X_Q$ coincides with the linear system $\Pi_Q.$
Hence to consider the complete intersection $X_Q$ of $m+1$ quadrics 
in $\projsp^n$ is equivalent to consider the linear system of quadrics defining $X_Q.$
In other words, the complete intersection $X_Q$
defines a unique $\GenLin_{m+1}(k) \times (k^\times I_{n+1})$-orbit in
$\Qdrcsp.$ 
Moreover, the projective equivalence class of $X_Q$ defines
a unique $\GenLin_{m+1}(k) \times \GenLin_{n+1} (k)$-orbit in $\Qdrcsp.$

\begin{rmk}
	Note that an element $M \in \RedLoci$ does not necessarily 
	define a complete intersection of 
	$m+1$ quadrics. An extreme example is given as follows.
	Let $Q_0, Q_1$ be two quadrics whose intersection is a smooth complete intersection.
	Then the discriminant polynomial of the pair $(Q_0, Q_1)$ is 
	a separable binary form (cf.\ \cite{Wang-Quadrics}).
	We put $Q' := (Q_0, Q_0, Q_1).$ 
	The corresponding triple of symmetric matrices $M_{Q'}$ is an element of
	$\RedLoci.$ But obviously it does not define a complete intersection of three quadrics.
\end{rmk}

\subsection{The projective automorphism groups of complete intersections of quadrics}
We fix an element $M_Q \in \NVLoci,$ and assume that
$M_Q$ defines a complete intersection $X_Q$ of $m+1$ quadrics in $\projsp^n.$
Let us write 
\[
	\phi_c(M_Q)= [(\mathcal{M}, \lambda, s)] \in \Triples.
\]
Here $\mathcal{M}$ is a coherent $\strshf_{\Pi_Q}$-module with
$\Supp(\mathcal{M}) = \Delta_Q.$
We put $L_0 = \Endomorphism{\Pi_Q}{\mathcal{M}}$
and use notations introduced in Section \ref{fibers} freely.
In particular, $L \subset L_0$ is the subspace of $L_0$ defined by the condition 
$\transpose l \circ \lambda = \lambda \circ l.$
We can define the norm map $\Norm$ as
\[
	\Norm \colon L_0 \longrightarrow L \quad ; \quad  l \mapsto \sigma(l)l.
\]
Since $\Norm(k^\times) \subset k^\times,$ the norm map induces the map
\[
	\overline{\Norm} \colon L_0^\times / k^\times \longrightarrow 
	(L \cap L_0^\times) / k^\times.
\]
This map is not a group homomorphism in general, but the kernel is always
a subgroup of $L_0^\times/k^\times.$

The projective automorphism group of $X_Q$ is defined by
\begin{align*}
	\Automorphism{\projsp^n}{X_Q} &:=
	\left\{g \in \Automorphism{}{\projsp^n} \; \relmid \; gX_Q = X_Q \right\}.
\end{align*}
We define another group
$\Automorphism{\Pi_Q}{\Delta_Q, \mathcal{M}, \lambda}$ as
\[
	\Automorphism{\Pi_Q}{\Delta_Q, \mathcal{M}, \lambda} :=
	\left\{ A \in \GenLin_{m+1}(k) \; \relmid \;
	\begin{aligned}
	&\nu_A(\Delta_Q) = \Delta_Q, \\
	&\Isomorphism{\Pi_Q}{(\nu_A^*\mathcal{M}, \nu_A^*\lambda),
	(\mathcal{M}, \lambda)} \neq \emptyset \ 
	\end{aligned}
	\right\} / (k^\times I_{m+1}),
\]
where we write
\[
	\Isomorphism{\Pi_Q}{(\mathcal{M}, \lambda), (\mathcal{M}', \lambda')}
	:= \left\{ \rho \colon \mathcal{M} \isomarrow \mathcal{M}' \; \relmid \; 
	\transpose \rho  \circ \lambda' \circ \rho = \lambda \right\}_.
\]

Now we can state our main result in this section.

\begin{thm}\label{ProjAutComm}
	Let $M_Q \in \NVLoci$ be an element defining 
	a complete intersection $X_Q$ of $m+1$ quadrics in $\projsp^n.$
	We write $\phi_c(M_Q) = [(\mathcal{M}, \lambda, s)] \in \Triples.$
	Then there exists a short exact sequence of the following form:
		\[
			\xymatrix{
				0 \ar[r] &
				\Kernel(\overline{\Norm}) \ar[r] &
				\Automorphism{\projsp^n}{X_Q} \ar[r] &
				\Automorphism{\Pi_Q}{\Delta_Q, \mathcal{M}, \lambda} \ar[r] &
				0.
			}
		\]
\end{thm}

Beauville proved this theorem when $m=2, n \ge 3,$ $k$ is algebraically closed and
$X_Q$ is a smooth complete intersection of three quadrics
(\cite[Proposition 6.19]{Beauville-Prym}).

The proof of Theorem \ref{ProjAutComm} is given in the next subsection.

\begin{rmk}
	The condition that $M_Q$ defines a complete intersection is stronger than
	the condition needed in the proof.
	We assume that the following conditions on $X_Q$ hold: $\Pi_Q \cong \projsp^m,$ 
	and any quadric $Q'$ containing $X_Q$ {\it as a subscheme}
	 is an element of $\Pi_Q.$
	Then we also have a short exact sequence as in Theorem \ref{ProjAutComm}.
\end{rmk}

\subsection{A description of the projective automorphism groups}\label{OPofAut}
In this subsection, we give an interpretation of the projective automorphism group
$\Automorphism{\projsp^n}{X_Q}$ of a 
complete intersection $X_Q$ of quadrics as an extension of two groups. 
Then we give a description of each group and prove Theorem \ref{ProjAutComm}.

Let $X_Q$ be a complete intersection of $m+1$ quadrics in $\projsp^n.$
We fix an $(m+1)$-tuple of quadrics 
\[
	Q=(Q_0, Q_1, \dots, Q_m)
\] 
defining $X_Q.$
Then we have
\[
	\Pi_Q := \langle Q_0, Q_1, \dots, Q_m \rangle \cong \projsp^m.
\]
We fix a quasi-isomorphism of 
complexes of coherent $\strshf_{\Pi_Q}$-modules
\[
	c \colon \omega_{\Pi_Q} \isomarrow \strshf_{\Pi_Q}(-m-1)[m].
\]
Assume that the corresponding $\GenLin_{m+1}(k) \times \GenLin_{n+1}(k)$-orbit of $M_Q$
is contained in $\NVLoci.$
We shall study the structure of the projective automorphism group 
$\Automorphism{\projsp^n}{X_Q}$ of $X_Q.$

Note that $\Automorphism{}{\projsp^n} \cong \ProjGL_{n+1}(k).$
Since $X_Q$ is a complete intersection of $m+1$ quadrics, 
if we take $P \in \GenLin_{n+1}(k)$ such that 
$\overline{P} \in \Automorphism{\projsp^n}{X_Q},$ the action of $P$ preserves
the linear system $\Pi_Q \subset \projsp(\Sym_2 k^{n+1}).$ 
Hence there exists a matrix $A \in \GenLin_{m+1}(k)$ such that
\[
	\transpose P M_Q P = M_Q \cdot A.
\]
Thus we have
\[
	\Automorphism{\projsp^n}{X_Q} \cong
	\left\{ P \in \GenLin_{n+1}(k) \; \relmid \; 
	\transpose P M_Q P = M_Q \cdot A 
	\quad (\exists A \in \GenLin_{m+1}(k)) \right\}/ (k^\times I_{n+1}).
\]

Let us define the groups $G_Q, F'_Q, E_Q, F_Q, P_Q, H_Q$ by
\begin{align*}
	G_Q &:= \left\{ (A, P) \in \GenLin_{m+1}(k) \times \GenLin_{n+1}(k) \; \relmid \; 
	\transpose P M_Q P = M_Q \cdot A \right\}, \\
	F'_Q &:= \left\{ (u,P) \in k^\times \times \GenLin_{n+1}(k) \; \relmid \; 
	\transpose P M_Q P = uM_Q \right\} \lhd G_Q, \\
	E_Q &:= \left\{ P \in \GenLin_{n+1}(k) \; \relmid \; 
	\transpose P M_Q P = M_Q \cdot A 
	\quad (\exists A \in \GenLin_{m+1}(k)) \right\}, \\
	F_Q &:= \left\{ P \in \GenLin_{n+1}(k) \; \relmid \; 
	\transpose P M_Q P = uM_Q 
	\quad (\exists u \in k^\times) \right\} \lhd E_Q, \\
	P_Q &:= \left\{(a^2I_{m+1}, aI_{n+1}) \relmid a \in k^\times \right\} 
	\subset G_Q, \\
	H_Q &:= E_Q / F_Q. 
\end{align*}
The group $F'_Q$ is a normal subgroup of $G_Q$ and
$F_Q$ is a normal subgroup of $E_Q.$
Obviously, $E_Q$ has a canonical surjection onto $\Automorphism{\projsp^n}
{X_Q}.$ 
Additionally, we find $G_Q$ is isomorphic to $E_Q.$

\begin{lem}\label{proj-isom}
	The second projection $\projection\colon G_Q \longrightarrow E_Q$ is an isomorphism.
	In particular, it induces an isomorphism $F'_Q \isomarrow F_Q.$
\end{lem}

\begin{prf}
	For each $P \in E_Q,$ a matrix $A \in \GenLin_{m+1}(k)$ with
	\[
		\transpose P M_Q P = M_Q \cdot A
	\]
	is determined by the action of $P$ on the $k$-vector space spanned by
	$M_{Q_0}, M_{Q_1}, \dots, M_{Q_m}.$
	Hence there exists a unique matrix $A$ satisfying this condition. \qed
\end{prf}

By this lemma, we have $H_Q \cong G_Q / F'_Q.$
The kernel of 
\[
	\projection_1 \colon G_Q \longrightarrow \GenLin_{m+1}(k)
\]
is obviously contained in $F'_Q,$ hence we also have
\[
	H_Q \cong \projection_1(G_Q) / \projection_1(F'_Q).
\]

\begin{lem}\label{kerproj}
	The kernel of the composite map
	\[
		G_Q \underset{\projection}{\isomarrow} E_Q \longrightarrow 
		\Automorphism{\projsp^n}{X_Q}
	\]
	coincides with $P_Q.$
\end{lem}

\begin{prf}
	If $(A, P)$ is an element in the kernel, 
	we have $P = aI_{n+1}$ for some $a \in k^\times.$
	Then we have
	\(
		a^2 M_Q = M_Q \cdot A.
	\)	
	Hence we have
	\(
		A = a^2 I_{m+1}.
	\)
\qed\end{prf}

Thus we have the following commutative diagram with exact rows and exact columns:
\[
	\xymatrix{
		& 1 \ar[d] & 1 \ar[d] & & \\
		& P_Q \ar@{=}[r] \ar[d]& P_Q \ar[d]& & \\
		1 \ar[r] & F_Q \ar[r] \ar[d]& E_Q = G_Q \ar[d]\ar[r] &
		 H_Q \ar[r] \ar@{=}[d] &1 \\
		1 \ar[r] & F_Q/P_Q \ar[r] \ar[d]& \Automorphism{\projsp^n}{X_Q}
		 \ar[d]\ar[r] & H_Q \ar[r] &1 \\
		 &1& 1.& &
	}
\]
In order to prove Theorem \ref{ProjAutComm},
it is enough to find an interpretation of the groups $F_Q/P_Q$ and $H_Q.$

\begin{lem}\label{small}
	We have
	\[
	F_Q \cong \left\{ l \in L_0^\times \; \relmid \; \sigma(l)l \in k^\times \right\}.
	\]
\end{lem}

\begin{prf}
	Take $P \in F_Q$ and $u \in k^\times$ satisfying $\transpose P M_Q P = uM_Q.$
	Since we have $\transpose P M_Q = u M_Q P^{-1},$ we have
	\(
		l := (P, uP^{-1}) \in L_0
	\) 
	by Proposition \ref{End1}.
	In particular, $\sigma(l) l = u \in k^\times.$ 
	Conversely, if $l =(P, P') \in L_0$ satisfies $\sigma(l)l = u\in k^\times,$
	we immediately see $P=uP'^{-1}.$
	Hence we have $\transpose P M_Q P = uM_Q.$ \qed
\end{prf}

\begin{prf}[Proof of Theorem \ref{ProjAutComm}]
By Lemma \ref{small}, we have
\[
	F_Q / P_Q \cong \Kernel (\overline{\Norm}) 
\]
because $\overline{\Norm} \colon L_0^\times / k^\times \longrightarrow 
L_0^\times / k^\times$ is defined by $\overline{\Norm} (l) = \sigma(l)l.$ 

Next we consider the group $H_Q.$
Let $(A, P) \in G_Q$ be an element. We have isomorphisms
\begin{align*}
	\nu_A &\colon \Pi_Q \isomarrow \Pi_Q, \\
	\rho &\colon \nu_A^* \mathcal{M} \isomarrow \mathcal{M}
\end{align*}
satisfying
\(
	\transpose \rho \circ \lambda \circ \rho = \nu_A^* \lambda.
\)
This shows $\nu_A \in \Automorphism{}{\Pi_Q}$ induces 
a projective automorphism of $\Delta_Q = \Supp(\mathcal{M}) \subset \Pi_Q.$
Hence the equivalence class $\overline{A} \in \ProjGL_{m+1}(k)$ of $A$ is 
an element of
\[
	\Automorphism{\Pi_Q}{\Delta_Q, \mathcal{M}, \lambda} :=
	\left\{ A \in \GenLin_{m+1}(k) \; \relmid \;
	\begin{aligned}
	&\nu_A(\Delta_Q) = \Delta_Q, \\
	&\Isomorphism{\Pi_Q}{(\nu_A^*\mathcal{M}, \nu_A^*\lambda),
	(\mathcal{M}, \lambda)} \neq \emptyset \ 
	\end{aligned}
	\right\} / (k^\times I_{m+1}),
\]
If $(A, P)$ is an element of $F'_Q,$ then we have $\nu_A = \identity_{\Pi_Q}$
as a projective automorphism.
Conversely, if $\nu_A = \identity_{\Pi_Q}$ then we can write
$A=aI_{m+1}$ for some $a \in k^\times.$ Hence $A \in \projection_1(F'_Q).$
Thus $H_Q = \projection_1(G_Q) / \projection_1(F'_Q)$ is a subgroup of
\(
	\Automorphism{\Pi_Q}{\Delta_Q, \mathcal{M}, \lambda}.
\)

If we take an element $A \in \GenLin_{m+1}(k)$ 
whose equivalence class $\overline{A} \in \ProjGL_{m+1}(k)$ is an element of 
$\Automorphism{\Pi_Q}{\Delta_Q, \mathcal{M}, \lambda},$
there exists an isomorphism
\[
	\rho \colon \nu_A^*\mathcal{M} \isomarrow \mathcal{M}
\]
satisfying \(
	\transpose \rho \circ \lambda \circ \rho = \nu_A^* \lambda.
\)
By a similar argument to the proof of Proposition \ref{fiber},
we have a matrix $P \in \GenLin_{n+1}(k)$ such that
\(
	\transpose P M_Q P = M_Q \cdot A.
\)
This shows $(A, P) \in G_Q$ and $\nu_A \in \projection_1(G_Q) / \projection_1(F'_Q).$
Thus the proof of Theorem \ref{ProjAutComm} is complete.\qed

\end{prf}

\end{document}